\documentclass[11pt,a4paper,twoside]{amsart}
\usepackage{amssymb,amsmath,amsthm,graphicx,color}
\theoremstyle{plain}
\newtheorem{theorem}{Theorem}[section]

\newtheorem{definition}[theorem]{Definition}
\newtheorem{lemma}[theorem]{Lemma}
\newtheorem{corollary}[theorem]{Corollary}
\newtheorem{proposition}[theorem]{Proposition}

\theoremstyle{remark}
\newtheorem{remark}[theorem]{Remark}

\usepackage{hyperref}

\allowdisplaybreaks[4]
\numberwithin{equation}{section}

\newcommand{\C}{\mathbb{C}}
\newcommand{\R}{\mathbb{R}}
\newcommand{\Z}{\mathbb{Z}}

\newcommand{\F}{\mathcal{F}}

\renewcommand{\Re}{\operatorname{Re}}
\newcommand{\I}{\infty}
\newcommand{\abs}[1]{\left\lvert #1\right\rvert}
\newcommand{\norm}[1]{\left\lVert #1\right\rVert}
\newcommand{\Lebn}[2]{\left\lVert #1 \right\rVert_{L^{#2}}}

\def\({\left(}
\def\){\right)}
\def\<{\left\langle}
\def\>{\right\rangle}
\def\le{\leqslant}
\def\ge{\geqslant}

\def\d{{\partial}}

\newcommand{\eps}{\varepsilon}

\newcommand{\rre}{\mathbb{R}}
\newcommand{\pt}{\partial}

\DeclareMathOperator{\sign}{sign}

\begin{document}
\title[On well-posedness of generalized KdV equation]
{On well-posedness of generalized\\
Korteweg-de Vries equation\\
in scale critical $\hat{L}^r$ space}
\author[S.Masaki and J.Segata]
{Satoshi Masaki  
and Jun-ichi Segata}
\address{Laboratory of Mathematics\\
Institute of Engineering\\
Hiroshima University\\
Higashihiroshima Hirhosima, 739-8527, Japan}
\email{masaki@amath.hiroshima-u.ac.jp}

\address{Mathematical Institute, Tohoku University\\
6-3, Aoba, Aramaki, Aoba-ku, Sendai 980-8578, Japan}
\email{segata@m.tohoku.ac.jp}
\maketitle

\begin{abstract}
The purpose of this paper is to study local and 
global well-posedness of initial value 
problem for generalized Korteweg-de Vries (gKdV) 
equation in $\hat{L}^r=\{f\in{{\mathcal S}}'(\rre)|
\norm{f}_{\hat{L}^r}=\|\hat{f}\|_{L^{r'}}<\infty\}$. 
We show (large data) local well-posedness, small data global 
well-posedness, and small data scattering for 
gKdV equation in the scale critical $\hat{L}^r$ space.
A key ingredient is a Stein-Tomas type inequality for the Airy 
equation, which generalizes usual Strichartz' estimates for $\hat{L}^r$-framework.
\end{abstract}

\section{Introduction}

We consider initial value problem for the 
generalized Korteweg-de Vries (gKdV) equation
\begin{eqnarray}
\left\{
\begin{array}{l}
\displaystyle{
\pt_tu+\pt_x^3u=\mu\pt_x(|u|^{\alpha-1}u),
\qquad t,x\in\rre,}\\
\displaystyle{u(0,x)=u_{0}(x),
\qquad\qquad\qquad\ \  x\in\rre,}
\end{array}
\right.
\label{gKdV}
\end{eqnarray}
where $u:\rre\times\rre\to\rre$ is an unknown function, 
$u_{0}:\rre\to\rre$ is a given function, 
and $\mu\in\rre\backslash\{0\}$ and $\alpha>1$ 
are constants. We call that (\ref{gKdV}) is defocusing if 
$\mu>0$ and focusing if $\mu<0$. 

The class of equations (\ref{gKdV}) arises in several 
fields of physics. Eq. (\ref{gKdV}) with $\alpha=2$ 
is notable Korteweg-de Vries equation which models long waves 
propagating in a channel \cite{KV}.  
Eq. (\ref{gKdV}) with $\alpha=3$ is also well known as the 
modified Korteweg-de Vries equation which 
describes a time evolution for the  curvature of 
certain types of helical space curves \cite{L}. 

The equation (\ref{gKdV}) has the following scale invariance: 
if $u(t,x)$ is a solution to (\ref{gKdV}), then
\[
	u_{\lambda}(t,x):=\lambda^{\frac{2}{\alpha-1}} u(\lambda^3 t, \lambda x)
\]
is also a solution to (\ref{gKdV}) with a initial data $u_{\lambda}(0,x)=
\lambda^{\frac{2}{\alpha-1}} u_{0}(\lambda x)$ for any $\lambda>0$.
In what follows, a Banach space for initial data is referred to as a \emph{scale critical space}
if its norm is invariant under $u_0(x) \mapsto 
\lambda^{\frac{2}{\alpha-1}} u_{0}(\lambda x)$.

The purpose of this paper is to study (large data) local 
well-posedness, small data global well-posedness 
and scattering for (\ref{gKdV}) in a scale critical space 
$\hat{L}^{(\alpha-1)/2}$. For $r \in[1,\I]$, the function space 
$\hat{L}^r$ is defined by 
\begin{eqnarray*}
	\hat{L}^r=\hat{L}^r(\rre):=\{f\in{{\mathcal S}}'(\rre)|
	\norm{f}_{\hat{L}^r} =\|\hat{f}\|_{L^{r'}}<\infty\},
\end{eqnarray*}
where $\hat{f}$ stands for Fourier transform of $f$ with respect to space variable
and $r'$ denotes the H\"older conjugate of $r$.
We use the conventions $1'=\I$ and $\I'=1$.
Our notion of well-poseness contains of
existence, uniqueness, and continuity of the data-to-solution map. 
We also consider persistent property of the solution, that is, the solution describes 
a continuous curve in the function space $X$ 
whenever $u_{0}\in X$. 

Local well-posedness of the initial 
value problem (\ref{gKdV}) in a scale subcritical Sobolev space 
$H^{s}(\rre)$, $s>s_\alpha:=1/2-2/(\alpha-1)$, has been studied by many authors \cite{B,G,Guo,Kato,KPV,KPVq,KN,MR1},
where $s_{\alpha}$, a scale critical exponent, is unique number such that $\dot{H}^{s_\alpha}$ becomes scale critical.  
A fundamental work on local well-posedness is due 
to Kenig-Ponce-Vega \cite{KPV}. They proved that (\ref{gKdV}) is locally 
well-posed in $H^{s}(\rre)$ with $s>3/4$ ($\alpha=2$, $s_2=-3/2$), 
$s\ge1/4$ ($\alpha=3$, $s_3=-1/2$), $s\ge1/12$ ($\alpha=4$, $s_4=-1/6$) and 
$s\ge s_{\alpha}$ ($\alpha\ge5$).
Introducing  Fourier restriction norms, Bourgain \cite{B} 
obtained local (and global\footnote{
Since the equation (\ref{gKdV}) preserves $L^{2}$ norm 
of solution in $t$, local well-posedness in $L^{2}$ yields 
global well-posedness in $L^{2}$ if $\alpha <5$.})
well-posedness of the KdV equation (i.e., (\ref{gKdV}) with $\alpha=2$) in $L^{2}(\rre)$. 
In \cite{KPVq}, Kenig-Ponce-Vega improved the previous results for 
the KdV equation to $H^{s}(\rre)$ with 
$s>-3/4$. Further, Guo \cite{Guo} and Kishimoto 
\cite{KN} extended Kenig-Ponce-Vega's result  
in $H^{-3/4}(\rre)$ (See also Buckmaster-Koch \cite{BK} 
on the existence of weak solution to the KdV equation at $H^{-1}$.). 
Gr\"{u}nrock \cite{G} has shown 
 local well-posedness of the quartic KdV equation 
((\ref{gKdV}) with $\alpha=4$) in $H^{s}$ with $s>s_{4}$. 
Notice that all of the above results are based on 
contraction mapping principle for corresponding integral equation.
Hence, a data-solution map associated with (\ref{gKdV}) is Lipschitz continuous
\footnote{In fact, 
if the nonlinear term is analytic, then the 
data-solution map associated with (\ref{gKdV}) is analytic}.  

Concerning the well-posedness 
of (\ref{gKdV}) in the scale critical $\dot{H}^{s_\alpha}$ space, 
Kenig-Pone-Vega \cite{KPV} proved local well-posedness and 
 global well-posedness for small data 
in the scale critical space $\dot{H}^{s_{\alpha}}$  
when $\alpha\ge5$.
Since the scale critical exponent $s_{\alpha}$ is
negative in the mass-subcritical case $\alpha<5$, 
well-posedness of (\ref{gKdV}) in $\dot{H}^{s_{\alpha}}$ becomes rather a difficult problem. 
Tao \cite{T} proved  local well-posedness
and global well-posedness for small data for (\ref{gKdV})
with the quartic nonlinearity $\alpha=4$
in $\dot{H}^{s_{4}}$. 
Later on, the above results are extended to a homogeneous Besov space 
$\dot{B}^{s_\alpha}_{2,\infty}$  
by Koch-Marzuola \cite{KM} ($\alpha = 4$) and Strunk ($\alpha \ge 5$).
As far as we know, 
local well-posedness and small data global 
well-posedness of (\ref{gKdV}) in $\dot{H}^{s_{\alpha}}$ 
for the mass-subcritical case $\alpha <5$ was open
except for the case $\alpha=4$.

Local and global well-posedness 
for a class of nonlinear dispersive equation 
is currently being intensively investigated
also in the framework of $\hat{L}^r$ space.
For one dimensional nonlinear Schr\"{o}dinger equation, 
\begin{eqnarray}
\left\{
\begin{array}{l}
\displaystyle{
i\pt_tv-\pt_x^2v=\mu |v|^{\alpha-1}v,
\qquad t,x\in\rre,}\\
\displaystyle{v(0,x)=v_{0}(x),
\qquad\qquad\quad\   x\in\rre,}
\end{array}
\right.
\label{NLS}
\end{eqnarray}
where $\mu \in\rre\backslash\{0\}$, 
Gr\"{u}nrock \cite{G2} has shown local and global existence of solution to 
(\ref{NLS}) with $\alpha=3$ in $\hat{L}^r$. 
Hyakuna-Tsutsumi \cite{HT} extended 
Gr\"{u}nrock's result 
in $\hat{L}^r$ to all mass-subcritical case 
$1<\alpha<5$. 
Gr\"{u}nrock \cite{G1} and Gr\"{u}nrock-Vega \cite{GV} 
proved local and global existence result for
the modified KdV equation 
(i.e., (\ref{gKdV}) with $\alpha=3$) in $\hat{H}_{s}^{r}$, 
where $\hat{H}_{s}^{r}=\{f\in{{\mathcal S}}';
\|f\|_{\hat{H}_{s}^{r}}=\|(1+\xi^{2})^{s/2}
\hat{f}(\xi)\|_{L_{\xi}^{r'}}<\infty\}$. 
However, the above results are not in scale critical settings.

It would be interesting to compare the scale critical space 
$\hat{L}^{\frac{\alpha-1}2}$ with some other
scale critical spaces in view of symmetries.
Other than the scaling,
the $\hat{L}^{\frac{\alpha-1}2}$-norm is invariant under the following 
three group operations
\begin{enumerate}
\item Translation in physical space: $(T_a f)(x) = f(x-a)$, $a\in \R$,
\item Translation in Fourier space: $(P_\xi f)(x) = e^{-ix\xi} f(x)$, $\xi\in\R$,
\item Airy flow: $(\mathrm{Ai}(t) f)(x) = e^{-t\d_x^3}f(x)$, $t\in \R$.
\end{enumerate}
The critical Lebesgue space $L^{\frac{\alpha-1}2}$ is invariant under 
the former two symmetries but not under the Airy flow.
The critical Sobolev space $\dot{H}^{s_\alpha}$ 
(or homogeneous Triebel-Lizorkin and homogeneous Besov spaces $\dot{A}^{s_\alpha}_{2,q}$ ($1\le q\le\I $), more generally) is not invariant with respect to $P_\xi$ if $s_\alpha \neq0$.
The critical weighted Lebesgue space $\dot{H}^{0,-s_\alpha}:=L^2(\R, |x|^{-2s_\alpha}dx)$ is not
invariant with respect to $T_a$ and $\mathrm{Ai}(t)$.
Further, when $\alpha=5$ these four spaces coincide with $L^2$, which is
invariant under the above three symmetries.
Thus, among the above four critical spaces, $\hat{L}^{\frac{\alpha-1}2}$ possesses
the most rich symmetries, and, in some sense, $\hat{L}^{\frac{\alpha-1}2}$ is close to $L^2$ space.
Inclusion relations between these spaces are summarized in Appendix \ref{sec:embedding}.

\subsection{Local well-posedness}
Before we state our main results, we 
introduce several notation. 

\begin{definition}
Let $(s,r) \in \R \times [1,\I]$. 
A pair $(s,r)$ 
is said to be \emph{acceptable} if
$1/r\in [0, 3/4)$ and 
\[
	s \in 
	\begin{cases}
	[-\frac{1}{2r} , \frac2r ] & 0\le \frac1r \le \frac12 ,\\
	(\frac2r - \frac54 , \frac52 - \frac3r) & \frac12 < \frac1r < \frac34.
	\end{cases}
\]
\end{definition}

For an interval $I \subset \R$ and an acceptable pair $(s,r)$, 
we define a function space $X(I;s,r)$ of space-time functions with the following norm
\[
	\norm{f}_{X(I;s,r)} = \norm{ |D_x|^{s} f }_{L^{p(s,r)}_x(\R; L^{q(s,r)}_t(I))} ,
\]
where the exponents in the above norm are given by
\begin{equation}\label{def:pq}
	\frac2{p(s,r)} + \frac1{q(s,r) }=\frac1r, \quad
	-\frac1{p(s,r)} + \frac2{q(s,r) } = s,
\end{equation}
or equivalently, 
\[
	\begin{pmatrix}
	1/{p(s,r)}\\ 
	1/{q(s,r)}
	\end{pmatrix}
	=\begin{pmatrix}
	-1/5 & 2/5\\ 
	2/5 & 1/5
	\end{pmatrix}
	\begin{pmatrix}
	s\\ 
	1/r
	\end{pmatrix}.
\]
We refer $X(I;s,r)$ to as an 
\emph{$\hat{L}^r$-admissible space}.

Our main theorems are as follows.

\vskip2mm

\begin{theorem}[local well-posedness in
$\hat{L}^{\frac{\alpha-1}2}$]\label{thm:lwp}
For $21/5<\alpha<23/3$, the problem (\ref{gKdV}) is locally well-posed in
$\hat{L}^{\frac{\alpha-1}2}$.
Namely,
for any $u_{0}\in \hat{L}_{x}^{\frac{\alpha-1}{2}}(\rre)$,
there exists an interval $I=I(u_0)$ such that
a unique solution 
\begin{eqnarray}
u\in C(I;\hat{L}_{x}^{\frac{\alpha-1}{2}}(\rre)) \cap 
\bigcap_{(s,\frac{\alpha-1}2):acceptable}X(I;s,\frac{\alpha-1}2)
\label{sol}
\end{eqnarray}
to (\ref{gKdV}) exists. Furthermore, for any given subinterval $I'\subset I$, 
there exists a neighborhood $V$ of $u_{0}$ in $
\hat{L}_{x}^{\frac{\alpha-1}{2}}(\rre)$ such that 
the map $u_{0}\mapsto u$ from $V$ into the class defined 
by (\ref{sol}) with $I'$ instead of $I$ is Lipschitz 
continuous. 
\end{theorem}

\vskip2mm

\begin{remark}
Theorem \ref{thm:lwp} (and all results below) holds for more general nonlinearity 
of the form $\d_x G(u)$ with $G \in \mathrm{Lip}\alpha$.
For precise condition on $G$, see Remark \ref{rem:generalG}.
\end{remark}

\vskip2mm

The proof of Theorem \ref{thm:lwp} is based on a contraction argument,
with a help of a space-time estimate for the Airy equation in $\hat{L}^{r}$. 
A key ingredient is Stein-Tomas type inequality for the Airy 
equation, a special case of \cite[Corollary 3.6]{G1}:
\begin{equation}
\norm{|D_x|^{1/r} e^{-t\d_x^3} 
f}_{L^r_{t,x}(I\times\rre)} \le C 
\norm{f}_{\hat{L}^{r/3}},
\label{SRF}
\end{equation}
where $r\in (4,\I]$. 
This inequality is a generalization of a well-known 
Strichartz estimate
\begin{equation*}
\norm{|D_x|^{1/6} e^{-t\d_x^3} 
f}_{L^6_{t,x}(I\times\rre)} \le C 
\norm{f}_{L^{2}}.
\end{equation*}
Moreover, interpolations 
between the above Stein-Tomas type inequality (\ref{SRF})
and
Kenig-Ruiz estimate or Kato's local smoothing effect
give us the following generalized Strichartz' estimate for the Airy equation
in $\hat{L}^r$-framework (Proposition \ref{ho}):
If $(s,r)$ is an acceptable pair then there exists $C$ such that
\begin{equation}\label{eq:ghS}
	\norm{ e^{-t\d_x^3}f}_{X(\R;s,r)} \le C 
	\norm{f}_{\hat{L}^{r}}
\end{equation}
for $f\in \hat{L}^r$.
Furthermore, 
combining the homogeneous estimate and Christ-Kiselev 
lemma (Lemma \ref{Christ-Kiselev}), we also obtain 
a generalized version of inhomogeneous Strichartz' estimates. 
The estimate \eqref{SRF} can be regarded as a kind of restriction estimate of Fourier transform,
which goes back to Stein and Tomas \cite{F,To}
(for more information on the restriction theorem, see e.g. \cite{TVV}).
It is worth mentioning that the $\hat{L}^r$ spaces have naturally come out in this context.

\vskip2mm

We set $S(I;r) := X(I;0,r)$. 
The $S(I;r)$ norm is so-called \emph{scattering norm}.
It is understood that a key for
obtaining a closed estimate for the corresponding integral equation, from which 
local well-posedness immediately follows, is to bound the scattering norm $S(I; \frac{\alpha-1}2)$.
In the proof of Theorem \ref{thm:lwp}, 
the scattering norm is handled by means of the above generalized Strihcartz' estimate \eqref{eq:ghS}.
Notice that the pair $(0,\frac{\alpha-1}{2})$ is acceptable only if $\alpha>21/5$.
Our restriction $\alpha > 21/5$ comes from this fact.
For the upper bound on $\alpha$, see Remark \ref{rem:upperbound}, below.
Alternatively, Sobolev's embedding also yields a bound on the scattering norm,
provided $\alpha \ge 5$.
In such case, we obtain local well-posedness in $\dot{H}^{s_\alpha}$ as in \cite{KPV}
(see Remark \ref{rem:HsLWP}).


\vskip2mm

\subsection{Persistence of regularity}
We establish two persistence-of-regularity type results
for $\hat{L}^{\frac{\alpha-1}2}$-solutions given in Theorem \ref{thm:lwp}. 
More specifically, we consider persistence of $\hat{L}^{r}$-regularity for $r \neq \frac{\alpha-1}2$
and $\dot{H}^s$ regularity for $-1<s< \alpha$.
These results yield local well-posedness in other $\hat{L}^r$ like space
such as $\hat{L}^{r_1} \cap \hat{L}^{r_2}$, $r_1 \le \frac{\alpha-1}2 \le r_2$,
and $\dot{H}^s \cap \hat{L}^{\frac{\alpha-1}2}$.

\vskip2mm

\begin{theorem}
[persistence of $\hat{L}^r$-regularity]
\label{prop:reg1}
Assume $21/5<\alpha<23/3$. 
Let $u_{0}\in \hat{L}_{x}^{\frac{\alpha-1}2}(\rre)$ and let $u \in C(I;\hat{L}^{\frac{\alpha-1}2}(\rre))$
be a corresponding solution given in Theorem \ref{thm:lwp}.
If $u_0 \in \hat{L}^{\frac{\alpha_0-1}2}_x$ for some $21/5<\alpha_0<23/3$, $\alpha_0\neq \alpha$, 
then
\begin{eqnarray*}
u\in C(I;\hat{L}_{x}^{\frac{\alpha_0-1}2}(\rre)) \cap \bigcap_{(s,\frac{\alpha_0-1}2):acceptable} X(I;s,\frac{\alpha_0-1}2).
\end{eqnarray*}
\end{theorem}

\vskip2mm

\begin{theorem}
[persistence of $\dot{H}^s$-regularity]\label{prop:reg2}
Assume $21/5<\alpha<23/3$. 
Let $u_{0}\in \hat{L}_{x}^{\frac{\alpha-1}2}(\rre)$ and let $u \in C(I, \hat{L}^{\frac{\alpha-1}2}(\rre))$
be a corresponding solution given in Theorem \ref{thm:lwp}.
If $u_0 \in \dot{H}^{\sigma}_x(\rre)$ for some $-1 < \sigma <\alpha$, then
\begin{eqnarray*}
|D_x|^{\sigma} u\in C(I;L^2(\rre)) \cap \bigcap_{(s,2):acceptable} X(I;s,2).
\end{eqnarray*}
\end{theorem}

\vskip2mm

As a corollary, we obtain the following well-posedness results.

\begin{corollary} We have the following.

\vskip1mm
\noindent
(i) If $21/5<\alpha<23/3$ then \eqref{gKdV} is locally well-posed in $\hat{L}^{r_1} \cap \hat{L}^{r_2}$ as long as $8/5 < r_1 \le \frac{\alpha-1}2 \le r_2 < 10/3$.

\vskip1mm
\noindent
(ii) If $21/5 < \alpha < 5$ then \eqref{gKdV} is locally well-posed in $\dot{H}^{s_\alpha} \cap 
\hat{L}^{\frac{\alpha-1}2}$, where $s_\alpha = \frac12-\frac2{\alpha-1}$.
\end{corollary}

\vskip2mm

Since $\hat{L}^{\frac{\alpha-1}2} \subset \dot{H}^{s_\alpha}$ does not hold (see Lemma \ref{lem:ws}),
the second is weaker than well-posedness in $\dot{H}^{s_\alpha}$.

Here we remark that an $\hat{L}^{\frac{\alpha-1}2}$-solution has conserved quantities,
provided the solution has appropriate regularity.
More precisely,
when $u_0\in \hat{L}^{\frac{\alpha-1}2} \cap L^2$, a solution $u(t)$ has a conserved mass 
\begin{eqnarray*}
M[u(t)] := \norm{u(t)}_{L^2}^2.
\end{eqnarray*}
Similarly, if $u_0 \in \hat{L}^{\frac{\alpha-1}2} \cap \dot{H}^1$ then energy 
\begin{eqnarray*}
E[u(t)] := \frac12 \Lebn{\d_x u(t)}2^2 + \frac{\mu}{\alpha+1} \norm{u(t)}_{L^{\alpha+1}}^{\alpha+1}
\end{eqnarray*}
is  invariant.
\subsection{Blowup and scattering}

We next consider long time behavior of solutions given in Theorem \ref{thm:lwp}.
To this end, 
we give the definitions of blow up and scattering of (\ref{gKdV}) 
for the initial data $u_{0}\in\hat{L}_{x}^{r}$. 
Set 
\begin{eqnarray*}
T_{\mathrm{max}}:&=&\sup\{T>0;\exists u\in C([0,T];\hat{L}_{x}^{r}(\rre))
:\text{solution\ to\ (\ref{gKdV})}\},\\
T_{\mathrm{min}}:&=&\sup\{T>0;\exists u\in C([-T,0];\hat{L}_{x}^{r}(\rre))
:\text{solution\ to\ (\ref{gKdV})}\}.
\end{eqnarray*}
Denote the lifespan of $u(t)$ as $(-T_{\mathrm{min}},T_{\mathrm{max}})$.
We say a solution $u(t)$ blows up in finite time for positive (resp. negative) time direction
if $T_{\mathrm{max}}<+\I$ (resp. $T_{\mathrm{min}}<+\I$). 
We say a solution $u(t)$ scatters for positive time direction
if $T_{\mathrm{max}}=+\I$ and there exists a unique function 
$u_{+}\in\hat{L}_{x}^{r}$ such that 
\begin{eqnarray*}
\lim_{t\to+\infty}\|u(t)-e^{-t\pt_{x}^{3}}u_{+}
\|_{\hat{L}_{x}^r}=0,
\end{eqnarray*}
where $e^{-t\d_x^3}u_{+}$ is a 
solution to the Airy equation $\pt_{t}v+\pt_{x}^{3}v=0$ 
with a initial condition $v(0,x)=u_{+}$. 
The scattering of $u$ for negative time direction is defined by a similar fashion. 

Roughly speaking, a solution scatters if linear dispersion effect dominates
the nonlinear interaction.
A typical case is when the data (and the corresponding solution) is small.
Here, we state this small data scattering for (\ref{gKdV}).
\begin{theorem}[Small data scattering]\label{thm:sg} 
Let $21/5<\alpha<23/3$. There exists $\varepsilon_{0}>0$ 
such that if $u_{0}\in \hat{L}_{x}^{\frac{\alpha-1}{2}}(\rre)$ 
satisfies $\|u_{0}\|_{\hat{L}_{x}^{\frac{\alpha-1}{2}}}\le\varepsilon_{0}$, then the solution $u(t)$ to (\ref{gKdV}) given in Theorem \ref{thm:lwp}
is global in time and scatters for both time directions.
Moreover,
\begin{eqnarray*}
\|u\|_{L_{t}^{\infty}(\R; \hat{L}_{x}^{\frac{\alpha-1}{2}})}
+\|u\|_{S(\R;\frac{\alpha-1}2)}\le 2 \|u_{0}\|_{\hat{L}_{x}^{\frac{\alpha-1}{2}}}.
\end{eqnarray*} 
\end{theorem}

We now give criterion for blowup and scattering.

\begin{theorem}[Blowup criterion]\label{prop:bc}
Assume $21/5<\alpha<23/3$. Let $u_0\in \hat{L}^{\frac{\alpha-1}2}$ and
let $u(t)$ be a corresponding unique solution of \eqref{gKdV} given in
Theorem \ref{thm:lwp}. 
If $T_{\mathrm{max}}<\I$ then 
\[
	\| u \|_{S([0,T);\frac{\alpha-1}2)} \to \I
\]
as $T \uparrow T_{\mathrm{max}}$.
A similar statement is true for negative time direction.
\end{theorem}

\vskip2mm

\begin{theorem}[Scattering criterion]\label{prop:sc}
Assume $21/5<\alpha<23/3$. Let $u_0\in \hat{L}^{\frac{\alpha-1}2}$ and
let $u(t)$ be a corresponding unique solution of \eqref{gKdV} given in
Theorem \ref{thm:lwp}. 
The solution $u(t)$ scatters forward in time if and only if
$T_{\mathrm{max}}=+\I$ and $\| u \|_{S([0,\I);\frac{\alpha-1}2)} <\I$.
A similar statement is true for negative time direction.
\end{theorem}


\vskip2mm

Finally, we give a criteria for scattering in terms of the energy.
We note that if an $\hat{L}^{\frac{\alpha-1}2}$-solution $u(t)$ scatters (in $\hat{L}^{\frac{\alpha-1}2}$ sense) 
as $t\to\pm\I$ and if $u_0 \in  \hat{L}^{\frac{\alpha_0-1}{2}}$
(resp. if $u_0 \in \dot{H}^\sigma$) then 
 $u(t)$ scatters as $t\to\pm\I$ also in $\hat{L}^{\frac{\alpha_0-1}{2}}$ sense
(resp. $\dot{H}^\sigma$ sense).

\begin{theorem}\label{thm:negativeE}
Let $21/5<\alpha<23/3$.
If $u_0 \in \hat{L}^{\frac{\alpha-1}2} \cap H^1$ satisfies $u_0 \neq0$ and $E[u_0] \le 0$
then $u(t)$ does not scatter as $t\to \pm \I$.
\end{theorem}

\vskip2mm

The rest of the paper is organized as follows. 
In Section 2, we prove some linear space-time estimates 
for solutions to the Airy equation, in $\hat{L}^r$-framework. 
The generalized Stirchartz estimates are established in Propositions \ref{ho} and \ref{inho}.
Section 3 is devoted to several nonlinear estimates.
We also introduce several function spaces to work with in this section.
Then, in Section 4, we prove our theorems.
In Appendix A, we prove a fractional chain rule in space-time function space (Lemma \ref{ff}). 
Finally in Appendix B, we briefly collect some inclusion relation for $\hat{L}^r$.

The following notation will be used throughout this 
paper: $|D_x|^s=(-\pt_x^2)^{s/2}$ and 
$\langle D_x\rangle^s=(I-\pt_x^2)^{s/2}$ denote 
the Riesz and Bessel potentials of order $-s$, 
respectively. For $1\le p,q\le\infty$ and $I\subset\rre$, 
let us define a space-time norm
\begin{eqnarray*}
\|f\|_{L_t^qL_x^p(I)}&=&
\|\|f(t,\cdot)\|_{L_x^p(\rre)}\|_{L_t^q(I)},\\
\|f\|_{L_x^pL_t^q(I)}&=&
\|\|f(\cdot,x)\|_{L_t^q(I)}\|_{L_x^p(\rre)}.
\end{eqnarray*}

\section{Linear Estimates for Airy Equation}

In this section we consider the space-time estimates 
of solution to the Airy equation
\begin{eqnarray}
\left\{
\begin{array}{l}
\displaystyle{
\pt_tu+\pt_x^3u=F(t,x),
\qquad\  t\in I, x\in\rre,}\\
\displaystyle{u(0,x)=f(x),
\qquad\qquad\ \ x\in\rre,}
\end{array}
\right.
\label{eq:A}
\end{eqnarray}
where $I\subset\rre$ is an interval, 
$F:I\times\rre\to\rre$ and $f:\rre\to\rre$ are 
given functions. 

Let $\{e^{-t\d_x^3}\}_{t\in\rre}$ be an isometric isomorphism 
group in $\hat{L}^r$ defined by $e^{-t\d_x^3} = \F^{-1} e^{it \xi^3} \F $,
or more precisely by
\begin{eqnarray*}
(e^{-t\d_x^3}f)(x)
=\frac{1}{\sqrt{2\pi}}
\int_{-\infty}^{\infty}e^{ix\xi+it\xi^{3}}
\hat{f}(\xi)d\xi.
\end{eqnarray*}
Using the group, 
the solution to (\ref{eq:A}) can be written as 
\begin{eqnarray*}
u(t)=e^{-t\d_x^3}f+\int_{0}^{t}e^{-(t-t')\d_x^3}F(t')dt'.
\end{eqnarray*}

We first show a homogeneous estimates associated with (\ref{eq:A}).

\vskip2mm

\begin{proposition}\label{ho} 
Let $I$ be an interval. 
Let $(p,q)$ satisfy
\[
	0\le\frac{1}{p}<\frac14,\quad\quad 0 \le \frac1q < \frac12 - \frac1p.
\]
Then, for any $f\in\hat{L}^r$,
\begin{equation}\label{eq:mixed}
	\norm{|D_x|^s e^{-t\d_x^3} f}_{L^p_x L^q_t(I)} \le C\norm{f}_{\hat{L}^r},
\end{equation}
where
\[
	\frac1r = \frac2p + \frac1q,\quad s=-\frac1p+\frac2q.
\]
and positive constant $C$ depends 
only on $r$ and $s$. 
\end{proposition}

\vskip2mm

Figure 1 shows 
the range of $(p,q)$ satisfying the assumption of Proposition \ref{ho}, where
$A=(1/4,0)$, $B=(1/4,1/4)$, 
and $C=(0,1/2)$. 
The line segments $OA$ and $OC$ 
is included, but the other 
parts of border are excluded.

\vskip5mm

\begin{center}
\includegraphics{Airy1.eps}
\end{center}
\vskip3mm
\begin{center}
Figure 1
\end{center}

\vskip5mm

To prove Proposition \ref{ho}, we 
show three lemmas. The first one 
is a Stein-Tomas type estimate.

\vskip2mm

\begin{lemma}[Stein-Tomas type estimate]
\label{ST}
For any $r\in (4,\I]$, 
there exists a positive constant $C$ depending 
only on $r$ such that for any $f\in \hat{L}^{r/3}$
\begin{equation}\label{eq:ST}
\norm{|D_x|^{1/r} e^{-t\d_x^3} 
f}_{L^r_{t,x}(I)} \le C 
\norm{f}_{\hat{L}^{r/3}}.
\end{equation}
\end{lemma}

\vskip2mm
\noindent
{\bf Proof of Lemma \ref{ST}.} 
Although a more general version is proved in \cite[Corollary 3.6]{G1},
here we give a direct proof which is based on the fact that the exponents for space-variable and time-variable in the left hand side
coincide. 

It suffices to prove (\ref{eq:ST}) for the case 
$I=\rre$. 
For notational simplicity, we omit $\rre$. 
The case $r=\I$ follows from the 
Hausdorff-Young inequality. Let $r<\I$.
Squaring both sides, we may show that
\begin{eqnarray}
	\norm{||D_x|^{1/r} e^{-t\d_x^3} 
	f|^2}_{L^{r/2}_{t,x}} 
	\le C \norm{f}_{\hat{L}^{r/3}}^2.
	\label{z1}
\end{eqnarray}
The left hand side of (\ref{z1}) is equal to
\[
	\norm{\iint_{\rre^{2}} e^{ix(\xi-\eta) 
	+ it(\xi^3-\eta^3)} |\xi\eta|^{1/r} \hat{f}(\xi)\overline{\hat{f}(\eta)}\,d\xi d\eta 
	}_{L^{r/2}_{t,x}}.
\]
Changing variables by $a=\xi-\eta$ 
and $b=\xi^3-\eta^3$, we have
\begin{eqnarray*}
\lefteqn{\norm{||D_x|^{1/r} e^{-t\d_x^3} 
	f|^2}_{L^{r/2}_{t,x}}}\\
	&=&\norm{\iint_{\rre^{2}} e^{ixa + itb} 
|\xi\eta|^{1/r} \hat{f}(\xi)\overline{\hat{f}(\eta)}\frac1{3|\xi^2-\eta^2|} \,dadb }_{L^{r/2}_{t,x}}.
\end{eqnarray*}
We now use the Hausdorff-Young inequality to deduce that
\begin{eqnarray}
\lefteqn{\norm{||D_x|^{1/r} e^{-t\d_x^3} 
	f|^2}_{L^{r/2}_{t,x}}}\label{m2}\\
	&\le& C \norm{|\xi\eta|^{1/r} 
	\hat{f}(\xi)\overline{\hat{f}(\eta)}|\xi^2-\eta^2|^{-1} }_{L^{(r/2)'}_{a,b}}\nonumber\\
	&=& C \left\{ \iint_{\R^2} \frac{|\xi\eta|^{\frac1{r-2}} |\hat{f}(\xi)|^{\frac{r}{r-2}}|\hat{f}(\eta)|^{\frac{r}{r-2}}}
	{|\xi-\eta|^{\frac{2}{r-2}} |\xi + \eta|^{\frac{2}{r-2}}}
	\,d\xi d\eta \right\}^{1-\frac{2}r}.\nonumber
\end{eqnarray}
Notice that $r/2\ge2$.
We now split the integral region $\R^2$ into $\{\xi\eta\ge0\}$ and $\{\xi\eta<0\}$.
We only consider the first case, since the other can be treated essentially in the same way.
For $(\xi,\eta)$ with $\xi \eta \ge 0$, we have
$\xi\eta \le (\xi + \eta)^2/4$, and so
\begin{eqnarray}
\lefteqn{	\iint_{\xi\eta\ge0} \frac{|\xi\eta|^{\frac1{r-2}} |\hat{f}(\xi)|^{\frac{r}{r-2}}|\hat{f}(\eta)|^{\frac{r}{r-2}}}
	{|\xi-\eta|^{\frac{2}{r-2}} |\xi + \eta|^{\frac{2}{r-2}}}
	\,d\xi d\eta}\label{m3}\\
	&\le& 
	C\iint_{\xi\eta\ge0} \frac{|\hat{f}(\xi)|^{\frac{r}{r-2}}|\hat{f}(\eta)|^{\frac{r}{r-2}}}
	{|\xi-\eta|^{\frac{2}{r-2}}}
	\,d\xi d\eta.\nonumber
\end{eqnarray}
By the H\"older and the Hardy-Littlewood-Sobolev 
inequality, we have
\begin{eqnarray}
	\lefteqn{\iint_{\xi\eta\ge0} \frac{|\hat{f}(\xi)|^{\frac{r}{r-2}}|\hat{f}(\eta)|^{\frac{r}{r-2}}}
	{|\xi-\eta|^{\frac{2}{r-2}}}
	\,d\xi d\eta}\label{m4}\\
	&\le& \Lebn{|\hat{f}|^{\frac{r}{r-2}}}{\frac{r-2}{r-3}}
	\Lebn{(|\xi|^{-\frac2{r-2}} * |\hat{f}|^{\frac{r}{r-2}})}{r-2}\nonumber\\
	&\le& C \Lebn{\hat{f}}{\frac{r}{r-3}}^{\frac{2r}{r-2}} = C \norm{f}_{\hat{L}_{x}^{r/3}}^{\frac{2r}{r-2}}
	\nonumber
\end{eqnarray}
as long as $2/(r-2)<1$, that is, $r>4$.
Combining (\ref{m2}),(\ref{m3}) and (\ref{m4}), 
we obtain the result.
$\qquad\qed$

\vskip2mm

The second is Kenig-Ruiz type estimate \cite{KeRu}. 

\vskip2mm

\begin{lemma}[Kenig-Ruiz type estimate]\label{KR} 
There exists a universal constant $C$ such that
for any interval $I$ and any $f\in L^2$
\begin{equation}\label{eq:KR}
	\norm{|D_x|^{-\frac14} e^{-t\d_x^3} 
	f}_{L^4_x L^\I_t(I)}
	\le C\Lebn{f}2.
\end{equation}
\end{lemma}

\vskip2mm
\noindent
{\bf Proof of Lemma \ref{KR}.}
See \cite[Theorem 2.5]{KPV1}.
$\qquad\qed$

\vskip2mm

The last estimate is an $\hat{L}^q$ version of 
the Kato's local smoothing effect
\cite{Kato}. 

\vskip2mm

\begin{lemma}[Kato's smoothing effect]\label{KS}
For any $q\in[2,\I]$, there exists a positive 
constant $C$ depending only on $q$  
such that any interval $I$ and for any $f\in \hat{L}^q$
\begin{equation}\label{eq:K}
	\norm{|D_x|^{\frac2q} e^{-t\d_x^3} f}_{L^\I_x L^q_t(I)}
	\le C\norm{f}_{\hat{L}^q}.
\end{equation}
\end{lemma}

\vskip2mm
\noindent
{\bf Proof of Lemma \ref{KS}.} 
We show \eqref{eq:K} by slightly modifying 
the argument due to Kenig-Ponce-Vega 
\cite[Theorem 2.5]{KPV1}. 
We prove (\ref{eq:K}) for the case $I=\rre$ only.

The case $q=\infty$ 
is treated in Lemma \ref{ST}.
Hence, we may suppose that $q<\I$. A direct computation shows
\begin{eqnarray*}
|D_x|^{\frac2q} e^{-t\d_x^3} f
&=&
\frac{1}{\sqrt{2\pi}}\int_\R e^{ix\xi+it\xi^3} |\xi|^{\frac2q} 
\hat{f}(\xi)\, d\xi\\
&=&\frac{1}{3\sqrt{2\pi}}
\int_\R e^{ix\eta^{1/3} + it \eta} 
|\eta|^{\frac2{3q}}\eta^{-\frac23} \hat{f}(\eta^{\frac13}) \,
d\eta, 
\end{eqnarray*}
where we have used a change of variable $\eta=\xi^3$ to yield the last line.
Take $L^q_t$ norm and apply the Hausdorff-Young inequality to obtain
\[
	\norm{|D_x|^{\frac2q} e^{-t\d_x^3} f}_{L^q_t}
	\le C \norm{ e^{ix\eta^{1/3}} |\eta|^{\frac{2-q}{3q}} \hat{f}(\eta^{\frac13})}_{L^{q'}_\eta}
	\le C \norm{\hat{f}}_{L^{q'}}=C \norm{f}_{\hat{L}^q}.
\]
Since the right hand side is independent of $x$, we obtain (\ref{eq:K}). $\qquad\qed$

\vskip2mm
\noindent
{\bf Proof of Proposition \ref{ho}.}
Interpolating \eqref{eq:ST}, 
\eqref{eq:KR}, and \eqref{eq:K},
we obtain \eqref{eq:mixed}. $\qed$

\vskip2mm

Next we show an inhomogeneous estimates associated with (\ref{eq:A}).

\vskip2mm

\begin{proposition}\label{inho} 
Let $4/3<r<4$ and let 
$(p_{j},q_{j})$ ($j=1,2$) satisfy
\[
	0\le\frac{1}{p_{j}}<\frac14,\quad\quad 0 \le \frac{1}{q_{j}} < 
	\frac12 - \frac1p_{j}.
\]
Then, the inequalities 
\begin{equation}
\left\|
\int_0^te^{-(t-t')\pt_{x}^{3}}F(t')dt'
\right\|_{L_t^{\infty}(I;\hat{L}_x^{r})}
\le C_{1}
\||D_x|^{-s_{2}}
F\|_{L_x^{p_{2}'}L_{t}^{q_{2}'}(I)},
\label{k}
\end{equation}
and
\begin{equation}
\left\||D_x|^{s_{1}}
\int_0^te^{-(t-t')\pt_{x}^{3}}F(t')dt'
\right\|_{L_x^{p_{1}}L_{t}^{q_{1}}(I)}
\le C_{2}
\||D_x|^{-s_{2}}
F\|_{L_x^{p_{2}'}L_{t}^{q_{2}'}(I)} \label{l}
\end{equation}
hold for any $F$ satisfying $|D_x|^{-s_{2}}
F\in L_x^{p_{2}'}L_{t}^{q_{2}'}$, where
\[
	\frac1r = \frac{2}{p_{1}} + \frac{1}{q_{1}},
	\quad s_{1}=-\frac{1}{p_{1}}+\frac{2}{q_{1}}
\]
and
\[
	\frac{1}{r'} = \frac{2}{p_{2}} + \frac{1}{q_{2}},
	\quad s_{2}=-\frac{1}{p_{2}}+\frac{2}{q_{2}},
\]
where 
the constant $C_{1}$ depends on $r$, $s_{1}$ and $I$, and 
the constant $C_{2}$ depends on $r$, $s_{1}$, $s_{1}$ and $I$.
\end{proposition}

\vskip2mm

To prove Theorem \ref{inho}, we employ the 
following lemma which is essentially due to 
Christ-Kiselev \cite{CK}. The version of this 
lemma that we use is the one presented in 
Molinet-Ribaud \cite{MR}.

\vskip2mm

\begin{lemma}\label{Christ-Kiselev} Let $I\subset\rre$ 
be an interval and let 
$K:{{{\mathcal S}}}(I\times\rre)\to C(\rre^3)$. 
Assume that 
\begin{eqnarray*}
\left\|\int_{I}
K(t,t')F(t')dt'\right\|_{L_x^{p_1}L_t^{q_1}(I)}
\le C\|F\|_{L_x^{p_2}L_t^{q_2}(I)}
\end{eqnarray*}
for some $1\le p_1,p_2,q_1,q_2\le\infty$ 
with $\min(p_1,q_1)>\max(p_2,q_2)$. Then
\begin{eqnarray*}
\left\|\int_0^tK(t,t')F(t')dt'
\right\|_{L_x^{p_1}L_t^{q_1}(I)}
\le C\|F\|_{L_x^{p_2}L_t^{q_2}(I)}.
\end{eqnarray*}
Moreover the case $q_1=\infty$ and $p_2,q_2<\infty$ 
is allowed. 
\end{lemma}

\vskip2mm
\noindent
{\bf Proof of Lemma \ref{Christ-Kiselev}.}
See \cite[Lemma 2]{MR}. $\qed$

\vskip2mm
\noindent
{\bf Proof of Proposition \ref{inho}.}
We first prove the inequality (\ref{k}).
By the $\hat{L}^{r}$-unitarity of the group 
$\{e^{-t\pt_{x}^{3}}\}_{t\in\rre}$, 
the duality argument and Proposition \ref{ho}, 
we have 
\begin{eqnarray}
\lefteqn{\left\|\int_{0}^{t}e^{-(t-t')\pt_{x}^{3}}
F(t')dt'
\right\|_{\hat{L}_x^{r}}}
\label{t4}\\
&=&\left\|\int_{0}^{t}e^{t'\pt_{x}^{3}}
F(t')dt'
\right\|_{\hat{L}_x^{r}}
\nonumber\\
&=&\sup_{\|g\|_{\hat{L}_x^{r'}}=1}
\left[
\int_{-\infty}^{\infty}
\left\{
\int_{0}^{t}e^{t'\pt_{x}^{3}}
F(t',x)dt'\right\}
g(x)dx
\right]
\nonumber\\
&=&\sup_{\|g\|_{\hat{L}_x^{r'}}=1}
\left[
\int_{0}^{t}
\!\int_{-\infty}^{\infty}
|D_{x}|^{-s_{2}}F(t',x)
|D_{x}|^{s_{2}}e^{-t'\pt_{x}^{3}}g(x)dt'dx
\right]
\nonumber\\
&\le&\sup_{\|g\|_{\hat{L}_x^{r'}}=1}
\||D_x|^{-s_{2}}F\|_{L_x^{p_{2}'}L_t^{q_{2}'}(I)}
\||D_x|^{s_{2}}
e^{-t'\pt_{x}^{3}}g
\|_{L_x^{p_{2}}
L_t^{q_{2}}(I)}
\nonumber\\
&\le& C\sup_{\|g\|_{\hat{L}_x^{r'}}=1}
\||D_x|^{-s_{2}}F\|_{L_x^{p_{2}'}L_t^{q_{2}'}(I)}
\|g\|_{\hat{L}_x^{r'}}
\nonumber\\
&=&C\||D_x|^{-s_{2}}F\|_{L_x^{p_{2}'}L_t^{q_{2}'}(I)},
\nonumber
\end{eqnarray}
where the constant $C$ is independent of $t$. 
Hence we have (\ref{k}).

Next we prove the the inequality (\ref{l}). 
Since the case $r=2$ has already proved in 
\cite{KPV}, we consider the case where $r\neq2$. 
To prove (\ref{l}), it suffices to prove  
\begin{eqnarray}
\qquad\left\||D_x|^{s_{1}}
\int_{I}e^{-(t-t')\pt_{x}^{3}}
F(t')dt'\right\|_{
L_x^{p_{1}}L_t^{q_{1}}(I)}
\le C
\||D_x|^{-s_{2}}
F\|_{L_x^{p_{2}'}L_t^{q_{2}'}(I)}.\label{l1}
\end{eqnarray}
Indeed, since 
\begin{eqnarray*}
\min(p_{1},q_{1})=\left\{
\begin{array}{l}
\frac{r}{r-1}\quad(\frac43<r<2),
\\
r\qquad(2<r<4)
\end{array}
\right.
>
\max(p_{2}',q_{2}')=\left\{
\begin{array}{l}
r\qquad(\frac43<r<2),
\\
\frac{r}{r-1}\quad(2<r<4),
\end{array}
\right.
\end{eqnarray*}
we see that the combination of the Christ-Kiselev lemma 
(Lemma \ref{Christ-Kiselev}) with  
(\ref{l1}) implies (\ref{l}). Therefore 
we concentrate our attention on prove (\ref{l1}). By 
Proposition \ref{ho},  
\begin{eqnarray}
\lefteqn{\left\||D_x|^{s_{1}}
\int_{I}e^{-(t-t')\pt_{x}^{3}}
F(t')dt'\right\|_{L_x^{p_{1}}L_t^{q_{1}}(I)}}
\label{ay}\\
&=&
\left\||D_x|^{s_{1}}e^{-t\pt_{x}^{3}}
\int_{I}e^{t'\pt_{x}^{3}}
F(t')dt'\right\|_{
L_x^{p_{1}}
L_t^{q_{1}}(I)}
\nonumber\\
&\le& C\left\|\int_{I}e^{t'\pt_{x}^{3}}
F(t')dt'
\right\|_{\hat{L}_x^{r}}.
\nonumber
\end{eqnarray}
By the duality argument similar to (\ref{t4}), we obtain
\begin{eqnarray}
\left\|\int_{I}e^{t'\pt_{x}^{3}}
F(t')dt'
\right\|_{\hat{L}_x^{r}}
\le C\||D_x|^{-s_{2}}F
\|_{L_x^{p_{2}'}L_t^{q_{2}'}(I)}.
\label{by}\end{eqnarray}
Combining (\ref{ay}) and (\ref{by}), 
we obtain (\ref{l1}). $\qed$

\section{Nonlinear estimates}

In this section, we prove 
several nonlinear estimates which are used to prove 
main theorems. 
We introduce several function spaces.
Let us recall that a 
pair $(s,r)\in \R \times [1,\I]$ 
is said to be \emph{acceptable} if
$1/r\in [0, 3/4)$ and 
\[
	s \in 
	\begin{cases}
	[-\frac{1}{2r} , \frac2r ] & 0\le \frac1r \le \frac12 ,\\
	(\frac2r - \frac54 , \frac52 - \frac3r) & \frac12 < \frac1r < \frac34.
	\end{cases}
\]

\begin{definition}
Let $(s,r) \in \R \times [1,\I]$. 
A pair $(s,r)$ is said to be 
\emph{conjugate-acceptable} if
$(1-s,r')$ is acceptable, where $\frac1{r'}=1-\frac1r \in [0,1]$.
\end{definition}
\vskip3mm

\begin{center}
\includegraphics{gKdV_sr_range_both.eps}
\end{center}
\vskip3mm
\begin{center}
Figure 2
\end{center}

\vskip3mm

Figure 2 shows the ranges of acceptable pairs (quadrangle OABC)
and conjugate-acceptable pairs (quadrangle DEFG).
Here, $O=(0,0)$, $A=(1/2,-1/4)$, $B=(3/4,1/4)$, $C=(1/2,1)$,
$D=(1,1)$, $E=(1/2,5/4)$, $F=(1/4,3/4)$, and $G=(1/2,0)$.

For an interval $I \subset \R$ and a conjugate-acceptable pair 
$(s,r)$,
we define a function space $Y(I;s,r)$ by
\[
	\norm{f}_{Y(I;s,r)} = \norm{ |D_x|^{s} f }_{L^{\tilde{p}(s,r)}_x(\R; L^{\tilde{q}(s,r)}_t(I))} ,
\]
where the exponents are given by
\begin{equation}\label{def:tpq}
	\frac2{\tilde{p}(s,r)} + \frac1{\tilde{q}(s,r) }=2+\frac1r, \quad
	-\frac1{\tilde{p}(s,r)} + \frac2{\tilde{q}(s,r) } = s,
\end{equation}
or equivalently, 
\[
	\begin{pmatrix}
	1/{\tilde{p}(s,r)}\\ 
	1/{\tilde{q}(s,r)}
	\end{pmatrix}
	=\begin{pmatrix}
	-1/5 & 2/5\\ 
	2/5 & 1/5
	\end{pmatrix}
	\begin{pmatrix}
	s\\ 
	2+1/r
	\end{pmatrix}
	=
	\begin{pmatrix}
	1/{p(s,r)}\\ 
	1/{q(s,r)}
	\end{pmatrix}
	+
	\begin{pmatrix}
	4/5\\ 
	2/5
	\end{pmatrix}
	.
\]

With this terminology, Propositions \ref{ho} and \ref{inho} can be reformulated as follows:

\begin{proposition}\label{prop:ho_inho} 
Let $I$ be an interval. 

\vskip1mm
\noindent
(i) Let $(s,r)$ be an acceptable pair.
Then, there exists a positive constant $C$ depending only on $s$ and $r$ such that
\[
	\norm{e^{-t\d_x^3} f}_{L^\I(\R;\hat{L}^r)} + \norm{e^{-t\d_x^3} f}_{X(\R;s,r)} \le C_{s,r} \norm{f}_{\hat{L}^r}
\]
for any $f\in\hat{L}^r$.

\vskip1mm
\noindent
(ii) Let $(s_1,r)$ be an acceptable pair and let $(s_2,r)$ be a conjugate-acceptable pair.
Then, there exists a positive constant 
depending only on $s_i$ and $r$ such that for any $t_0 \in I \subset \R$ 
and any $F\in Y(I;s_2,r)$, 
\[
	\norm{ \int_{t_0}^t e^{-(t-t')\d_x^3 }\d_x F(t') dt'}_{L^\I_t(I;\hat{L}^r_x) \cap X(I;s_1,,r)}
	\le C\norm{F}_{Y(I;s_2,r)}.
\]
\end{proposition}

To handle $X(I;s,r)$ and $Y(I;s,r)$ spaces,
the following lemma is useful.
\begin{lemma}\label{lem:gHolder}
Let $1< p_i,q_i <\I$ and $s_i \in \R$ for $i=1,2$.
Let $p,q,s$ be 
\[
	\frac1p = \frac{\theta}{p_1}+\frac{1-\theta}{p_2}, \quad
	\frac1q = \frac{\theta}{q_1}+\frac{1-\theta}{q_2}, \quad
	s = \theta s_1+(1-\theta)s_2
\] 
for some $\theta \in (0,1)$.
Then, there exists a positive constant $C$ depending 
on $p_{1},p_{2},q_{1},q_{2},s_{1},s_{2}$ and $\theta$ such that
\[
	\norm{|D_x|^s f}_{L^p_x L^q_t}
	\le C \norm{|D_x|^{s_1} f}_{L^{p_1}_x L^{q_1}_t}^\theta \norm{|D_x|^{s_2} f}_{L^{p_2}_x L^{q_2}_t}^{1-\theta}
\]
holds 
for any $f$ such that $|D_x|^{s_1} f \in L^{p_1}_x L^{q_1}_t$ and $|D_x|^{s_2} f \in L^{p_2}_x L^{q_2}_t$.
\end{lemma}

\vskip2mm
\noindent
{\bf Proof of Lemma \ref{lem:gHolder}.}
For $z\in \C$,
define an operator $T_z = |D_x|^{zs_1+ (1-z)s_2 }$.
Let $g(t)$ and $h(x)$ be $\R$-valued simple functions and $G_z(t)$ and $H_z(x)$ be
extensions of these functions defined by
\[
	G_z (t) :=
	|g(t)|^{\frac{1-(z/{q_1}+(1-z)/q_2)}{1-1/q}} \sign g(t) 
\]
and
\[
	H_z (x) :=
	|h(x)|^{\frac{1-(z/{p_1}+(1-z)/p_2)}{1-1/p}} \sign h(x) ,
\]
respectively,
for $z\in\C$ with $0 \le \Re z \le 1$.
Put 
\[
	\Psi(z):=\iint_{\R^2} T_z f(t,x) {G_z(t)H_z(x)}dt dx.
\]
By density and duality, it suffices to show
\begin{equation}\label{eq:gHolder_pf1}
	\abs{ \Psi(\theta)} 
	\le C \norm{|D_x|^{s_1} f}_{L^{p_1}_x L^{q_1}_t}^\theta \norm{|D_x|^{s_2} f}_{L^{p_2}_x L^{q_2}_t}^{1-\theta}
\end{equation}
for any $f\in \mathcal{S}(\R^2)$ with compact Fourier support
and any simple functions $g(t)$ and $h(x)$ such that $\norm{g}_{L^{q'}_t} = \norm{h}_{L^{p'}_x}=1$. 

Let us prove \eqref{eq:gHolder_pf1}.
It is easy to see that $\Psi(z)$ is analytic in $0< \Re z < 1$ and continuous in $0\le \Re z \le 1$.
By a variant of multiplier theorem by Fernandez \cite[Theorem 6.4]{Fz}, 
we see that
$|D_x|^{it}$ is a bounded operator in $L^{p_1}_x L^{q_1}_t$ with norm $C(1+|t|)$.
Therefore, for any $y\in\R$,
\begin{align}
	|\Psi(1+iy)| &{}\le \norm{|D_x|^{iy(s_1-s_2)} (|D_x|^{s_1} f)}_{L^{p_1}_x L^{q_1}_t}
	\norm{G_{1+iy}H_{1+iy}}_{L^{p_1'}_x L^{q_1'}_t} 
	\label{w1}\\
	&{} \le C(1+|y(s_1-s_2)|) \norm{|D_x|^{s_1} f}_{L^{p_1}_x L^{q_1}_t}\norm{g}_{L^{q'}_t} \norm{h}_{L^{p'}_x}
	\nonumber\\
	&{} \le C(1+|y(s_1-s_2)|) \norm{|D_x|^{s_1} f}_{L^{p_1}_x L^{q_1}_t}.
	\nonumber
\end{align}
The same argument yields 
\begin{eqnarray}
|\Psi(iy)| \le C(1+|y(s_1-s_2)|) \norm{|D_x|^{s_2} f}_{L^{p_2}_x L^{q_2}_t}.
\label{w2}
\end{eqnarray}
From (\ref{w1}), (\ref{w2}) and Hirschmann's Lemma \cite{H}, 
we obtain \eqref{eq:gHolder_pf1} (see also \cite{S}). 
$\qed$

\subsection{Estimates on nonlinearity}
In this subsection, we establish an estimate on nonlinearity.
For this, we introduce a Lipschitz $\mu$ norm ($\mu>0$) as follows.
Write $\mu = N+\beta$ with $N \in \Z$ and $\beta \in (0,1]$.
For a function $G: \C \to \C$, we define
\[
	\norm{G}_{\mathrm{Lip} \mu} := 
	\sum_{j=0}^N \sup_{z \in \R\setminus \{0\}} \frac{|G^{(j)}(z)|}{|z|^{\mu-j}}+
	\sup_{x\neq y} \frac{| G^{(N)}(x)-G^{(N)}(y) |}{|x-y|^{\beta}}.
\]
where $G^{(j)}$ is $j$-th derivative of $G$. 
We say $G\in\mathrm{Lip} \mu$ if $G\in C^{N}(\R)$ and $\norm{G}_{\mathrm{Lip} \mu}<\infty$. 

\vskip2mm

The main estimates of this subsection is as follows:

\vskip2mm

\begin{lemma}\label{lem:nlest}
Suppose that $G(z) \in \mathrm{Lip} \alpha $ for some
$21/5<\alpha < 23/3$.
Let $(s,r)$ be a pair which is acceptable and conjugate-acceptable.
Then, the following two assertions hold:

\vskip1mm
\noindent
(i) If $u \in S(I;\frac{\alpha-1}2) \cap X(I;s,r)$ then $G(u) \in Y(I;s,r)$.
Moreover, there exists  a constant $C$ such that
\begin{align*}
\|G(u) \|_{Y(I;s, r)}  \le
C\norm{u}_{S(I;\frac{\alpha-1}2)}^{\alpha-1} \|u\|_{X(I;s,r)}
\end{align*}
for any $u \in S(I;\frac{\alpha-1}2)  \cap X(I;s,r)$. 

\vskip1mm
\noindent
(ii) There exists  a constant $C$ such that
\begin{align*}
 & \|G(u)-
G(v)\|_{Y(I;s, r)} \\
 &{} \le
 C (\|u\|_{X(I;s,r)} + \|v\|_{X(I;s,r)})\\
 &{} \quad \quad \times
(\norm{u}_{S(I;\frac{\alpha-1}2)} + \norm{v}_{S(I;\frac{\alpha-1}2)})^{\alpha-2} \norm{u-v}_{S(I;\frac{\alpha-1}2)} \\
&\quad +
 C(\norm{u}_{S(I;\frac{\alpha-1}2)} + \norm{v}_{S(I;\frac{\alpha-1}2)})^{\alpha-1}
\|u-v\|_{X(I; s ,r)}
\end{align*}
for any $u,v \in S(I;\frac{\alpha-1}2) \cap X(I;s,r)$. 
\end{lemma}

\begin{remark}\label{rem:generalG}
It is easy to see that $|z|^{\alpha-1} z \in \mathrm{Lip} \alpha$.
The validity of the above lemma is all assumption on the nonlinearity that we need.
Hence, the all results of this article hold for an equation with generalized nonlinearity $\pt_tu+\pt_x^3u= \pt_x (G(u))$, provided $G(z) \in \mathrm{Lip} \alpha$.
\end{remark}

To prove the above lemma, we recall the following two lemmas.

\begin{lemma}\label{pp} 
Let $I$ be an interval. 
Assume that $s \ge 0$. 
Let $p,q, p_{i},q_{i}, \in (1,\infty)$ ($i=1,2,3,4$).
Then, we have
\begin{multline*}
\||D_{x}|^{s}(fg)\|_{L_{x}^{p}L_{t}^{q}(I)} \le \\
C(\||D_{x}|^{s}f\|_{L_{x}^{p_{1}}L_{t}^{q_{1}}(I)}
\|g\|_{L_{x}^{p_{2}}L_{t}^{q_{2}}(I)}
+ \|f\|_{L_{x}^{p_{3}}L_{t}^{q_{3}}(I)}
\||D_{x}|^{s} g\|_{L_{x}^{p_{4}}L_{t}^{q_{4}}(I)})
\end{multline*}
provided that 
\begin{eqnarray*}
\frac{1}{p}=\frac{1}{p_{1}}+\frac{1}{p_{2}}=\frac{1}{p_{3}}+\frac{1}{p_{4}},
\quad\frac{1}{q}=\frac{1}{q_{1}}+\frac{1}{q_{2}}=\frac{1}{q_{3}}+\frac{1}{q_{4}},
\end{eqnarray*}
where the constant $C$ is independent of 
$I$ and $f$. 
\end{lemma}

\vskip2mm
\noindent
{\bf Proof of Lemma \ref{pp}.} 
If $s \in \Z$ then (classical) Leibniz' rule, H\"older's inequality, and Lemma \ref{lem:gHolder}
give us the result.
By a similar argument, it suffices to consider the case $0<s<1$ to handle the general case. 
However, that case follows from \cite[Theorem A.8]{KPV} and Lemma \ref{lem:gHolder}. $\qed$

\vskip2mm

\begin{lemma}\label{ff}
Suppose that $\mu>1$ and $s\in(0,\mu)$. Let $G\in\mathrm{Lip} \mu$. 
If $p,p_1,p_2,q,q_1,q_2\in(1,\I)$ satisfies
\[
	\frac1p = \frac{\mu-1}{p_1} + \frac1{p_2},
	\qquad \frac1q = \frac{\mu-1}{q_1} + \frac1{q_2},
\]
then there exists a positive constant $C$ depending on $\mu,s,p_1,p_2,q_1,q_2$ 
and $I$ 
such that 
\[
	\norm{|D_{x}|^{s}G(f)}_{L_{x}^{p}L_{t}^{q}(I)}
	\le C \norm{G}_{\mathrm{Lip} \mu}
	\norm{f}_{L_{x}^{p_{1}}L_{t}^{q_{1}}(I)}^{\mu-1} 
	\norm{|D_{x}|^{s}f}_{L_{x}^{p_{2}}L_{t}^{q_{2}}(I)}
\]
holds for any $f$ satisfying $f\in L_{x}^{p_{1}}L_{t}^{q_{1}}(I)$ and 
$|D_{x}|^{s}f\in L_{x}^{p_{2}}L_{t}^{q_{2}}(I)$.
\end{lemma}

Although Lemma \ref{ff} is essentially the same as \cite[Theorem A.6]{KPV}, 
we give the proof of this lemma in Appendix A for self-containedness and in order to 
clarify the necessity of the assumption $G \in \mathrm{Lip} \mu$.

\vskip2mm

\noindent
{\bf Proof of Lemma \ref{lem:nlest}.}
We prove the second assertion since the first immediately follows from the second
by letting $v=0$.
For simplicity, we name
$S=S(I;\frac{\alpha-1}2)$, $L=X(I;s,r)$, and $N=Y(I;s, r)$.

Let us write 
\[
G(u)-G(v)
= (u-v) \int_0^1 G'(\theta u + (1-\theta)v ) d\theta .
\]
Lemma \ref{pp} implies that
\begin{eqnarray*}
\lefteqn{\| G(u)-G(v) \|_{N}}\\
&\le& C\|u-v\|_{S}
\int_0^1 \||D_{x}|^{s}\{G'(\theta u + (1-\theta)v)\}
\|_{L_{x}^{p_1}L_{t}^{q_1}} d\theta
\\
& &+C\|u-v\|_{L}
\int_0^1 \|\{G'(\theta u + (1-\theta)v)\}
\|_{L_{x}^{p_2}L_{t}^{q_2}} d\theta\\
&=:&I_{1}+I_{2},
\end{eqnarray*}
where
\begin{align*}
	\begin{pmatrix}
	1/{p_1}\\ 
	1/{q_1}
	\end{pmatrix}
	&{}= 
	\begin{pmatrix}
	1/{\tilde{p}(s,r)}\\ 
	1/{\tilde{q}(s,r)}
	\end{pmatrix}-
	\begin{pmatrix}
	1/{p(0, \frac{\alpha-1}2)}\\ 
	1/{q(0, \frac{\alpha-1}2)}
	\end{pmatrix}
	\\
	&{}= (\alpha-2)
	\begin{pmatrix}
	1/{p(0, \frac{\alpha-1}2)} \\
	1/{q(0,\frac{\alpha-1}2)}
	\end{pmatrix}
	+ 
	\begin{pmatrix}
	1/{p( s, r)} \\
	1/{q( s ,r )}
	\end{pmatrix}
\end{align*}
and
\begin{align*}
	\begin{pmatrix}
	1/{p_2}\\ 
	1/{q_2}
	\end{pmatrix}
	&{}= \begin{pmatrix}
	1/{\tilde{p}(s,r)}\\ 
	1/{\tilde{q}(s,r)}
	\end{pmatrix}-
	\begin{pmatrix}
	1/{p( s, r)}\\ 
	1/{q(s, r)}
	\end{pmatrix} \\
	&{}= (\alpha-1)
	\begin{pmatrix}
	1/{p(0, \frac{\alpha-1}2)} \\
	1/{q(0,\frac{\alpha-1}2)}
	\end{pmatrix}.
\end{align*}
It is easy to see that $\norm{G'}_{\mathrm{Lip}(\alpha-1)}\le \norm{G}_{\mathrm{Lip}\alpha} <+\I$.
By definition of $\norm{\cdot}_{\mathrm{Lip}(\alpha-1)}$,
we estimate $I_2$ as
\begin{eqnarray*}
I_{2}
&\le&C\|u-v\|_{L} \norm{G'}_{\mathrm{Lip}(\alpha-1)} \int_0^1 \| |\theta u+(1-\theta)v|^{\alpha-1}
\|_{L_{x}^{p_2}L_{t}^{q_2}} d\theta \\
&\le&C\|u-v\|_{L} \int_0^1 (\|u\|_S+\|v\|_S)^{\alpha-1} d\theta \\
&\le&
C
(\|u\|_{S}
+\|v\|_{S})^{\alpha -1}
\|u-v\|_{L}.
\end{eqnarray*}
On the other hand,
we see from Lemma \ref{ff} that
\begin{multline*}
	\||D_{x}|^{s}\{G'(\theta u + (1-\theta)v)\}
\|_{L_{x}^{p_1}L_{t}^{q_1}} \\
	\le C\norm{G'}_{\mathrm{Lip}(\alpha-1)}
	\norm{\theta u+ (1-\theta) v}_{S}^{\alpha-2} \norm{\theta u+ (1-\theta) v}_{L}
\end{multline*}
for any $\theta \in (0,1)$. Hence, we find the following estimate on $I_1$;
\[
I_{1}
\le C\|u-v\|_{S}\norm{G'}_{\mathrm{Lip}(\alpha-1)}
	(\norm{u}_S + \norm{v}_{S})^{\alpha-2} (\norm{u}_L+\norm{v}_{L}).
\]
Collecting the above inequalities, we obtain the result.
$\qed$


\section{Proof of main theorems} 

In this section, we prove the main theorems. 
Recall the notation $S(I;r)=X(I;0,r)$.
Now, take a number $s_L(\alpha)$ so that 
a pair 
$(s_L(\alpha),\frac{\alpha-1}2)$ is acceptable and conjugate-acceptable.
We denote $L(I;\frac{\alpha-1}2)=X(I;s_L(\alpha),\frac{\alpha-1}2)$ and
$N(I;\frac{\alpha-1}2)=Y(I; s_L(\alpha),\frac{\alpha-1}2)$.
\begin{remark}\label{rem:upperbound}
If  $27/7 < \alpha < 23/3$ then $s_L(\alpha)$ with the above property exists.
Indeed, $s_L(\alpha) = \frac34 -\frac1{\alpha-1}$ works.
Our upper bound on $\alpha$ comes from this point.
\end{remark}

%
%

\subsection{Local well-posedness in a scale-critical space}

Let us prove Theorem \ref{thm:lwp}. To prove this theorem, 
we show the following lemma. 

\begin{lemma}\label{lem:lwp}
Assume $21/5<\alpha<23/3$
and $u_{0}\in \hat{L}_{x}^{\frac{\alpha-1}{2}}$.
Let $t_0 \in \R$ and $I$ be an interval with $t_{0}\in I$.
Then, there exists a universal constant $\delta>0$ such that,
if  a tempered distribution $u_0$ and an interval $I\ni t_0$ satisfy
\[
	\eps=\eps(I;u_0,t_0):=\norm{e^{-(t-t_0)\d_x^{3}} u_0}_{S(I;\frac{\alpha-1}2)}
	+ \norm{e^{-(t-t_0)\d_x^{3}}u_0}_{L(I,\frac{\alpha-1}2)}
	\le \delta,
\]
then there exists a unique solution $u\in C(I;\hat{L}_{x}^{\frac{\alpha-1}{2}})$ 
to the following initial value problem
\begin{eqnarray*}
\left\{
\begin{array}{l}
\displaystyle{
\pt_tu+\pt_x^3u=\mu\pt_x(|u|^{\alpha-1}u),
\qquad t,x\in\rre,}\\
\displaystyle{u(t_0,x)=u_{0}(x),
\qquad\qquad\qquad\ x\in\rre}
\end{array}
\right.
\end{eqnarray*}
(in the sense of corresponding integral equation) and satisfies
\[
	\norm{u}_{S(I;\frac{\alpha-1}2)}+ \norm{u}_{L(I;\frac{\alpha-1}2)} \le 2\eps.
\]
If $u_0 \in \hat{L}^{\frac{\alpha-1}2}$, in addition, then
$$\norm{u}_{L^\I(I;\hat{L}^{\frac{\alpha-1}2})} \le \norm{u_0}_{\hat{L}^{\frac{\alpha-1}2}} +
C\eps^\alpha$$ holds for some constant $C>0$
and $u$ belongs to all $\hat{L}^{\frac{\alpha-1}2}$-admissible space $X(I;s,\frac{\alpha-1}2)$.
\end{lemma}

\vskip2mm
\noindent
{\bf Proof of Lemma \ref{lem:lwp}.}
For $R>0$, define a complete metric space 
\begin{align*}
Z_{R}
={}&
\left\{ u\in L\(I;\frac{\alpha-1}2\) \cap S\(I;\frac{\alpha-1}2\)
;
\|u\|_{Z}\le R
\right\} ,\\
\|u\|_{Z}:={}&
\|u\|_{L(I;\frac{\alpha-1}2)}+\|u\|_{S(I;\frac{\alpha-1}2)}, \quad
d_Z(u,v):=
\|u-v\|_Z.
\end{align*}
For given tempered distribution $u_{0}$ with $e^{-(t-t_0)\d_x^3}u_0 \in Z_\delta$
and $v\in Z_{R}$, we denote 
\begin{eqnarray*}
\Phi(v)(t)
:=e^{-(t-t_0)\pt_{x}^{3}}u_{0}
+\mu\int_{t_0}^{t}
e^{-(t-t')\pt_{x}^{3}}\pt_x(|v|^{\alpha-1}v)
(t')dt'.
\end{eqnarray*}
We show that there exist $\delta>0$ such that
$\Phi:Z_{2\eps}\to Z_{2\eps}$ is a contraction map
for any $0<\eps\le \delta$.
 
To this end, we prove that there exist constants 
$C_{1},C_{2}>0$ such that for any $u,v\in Z_{R}$, 
\begin{eqnarray}
\|\Phi(u)\|_{Z}
&\le& \|e^{-(t-t_0)\d_x^3} u_{0}\|_{Z}+C_{1}R^{\alpha},\label{eq:lwp_b1}\\
d_Z(\Phi(u),\Phi(v))
&\le&C_{2}R^{\alpha-1} d_Z(u,v).\label{eq:lwp_b2}
\end{eqnarray}
Let $u\in Z_R$. 
We infer from Proposition \ref{prop:ho_inho} (ii) that
\begin{eqnarray*}
\|\Phi(u)\|_{Z}
\le \|e^{-t\d_x^3}u_0\|_Z +C
\| |u|^{\alpha-1}u\|_{N(I;\frac{\alpha-1}2)}.
\end{eqnarray*}
We then apply Lemma \ref{lem:nlest} (i) with $r=\frac{\alpha-1}2$ and
 $s=s_L(\alpha)$ to obtain \eqref{eq:lwp_b1}.
A similar argument shows \eqref{eq:lwp_b2}.
We just employ Lemma \ref{lem:nlest} (ii) instead.

Now let us choose $\delta>0$ so that
\begin{equation}\label{eq:lwp_pf1}
\quad C_{1}(2\delta)^{\alpha-1}\le\frac12, \qquad
C_{2}(2\delta)^{\alpha-1}\le\frac12, 
\end{equation}
Then, we conclude from \eqref{eq:lwp_b1}, \eqref{eq:lwp_b2}, and the smallness assumption
that $\Phi$ is a contraction 
map on $Z_{2\eps}$. Therefore, the Banach fixed point 
theorem ensures that there exists a unique solution 
$u\in Z_{2\eps}$ to \eqref{gKdV}. 

We now suppose that $u_0 \in \hat{L}^{\frac{\alpha-1}2}$.
By means of
Proposition \ref{prop:ho_inho}, we have
\[
	\norm{u}_{L^\I(I, \hat{L}^{\frac{\alpha-1}2})} 
	\le \norm{u_0}_{\hat{L}^{\frac{\alpha-1}2}} + C\eps^\alpha
\]
as in \eqref{eq:lwp_b1}. 
The same argument shows $u\in X(I;s,\frac{\alpha-1}{2})$ for any $s$ such that
$(s,\frac{\alpha-1}2)$ is acceptable.
$\qed$

\vskip2mm
\noindent
{\bf Proof of Theorem \ref{thm:lwp}.}
By Lemma \ref{lem:lwp}, we obtain a unique solution $$u
\in L^\I_t([-T,T];\hat{L}^{\frac{\alpha-1}2}_x) \cap S([-T,T];\frac{\alpha-1}2) \cap L([-T,T];\frac{\alpha-1}2)$$ 
for small $T=T(u_0)>0$.
We repeat the above argument to extend the solution,
and then obtain a solution which has a maximal lifespan. 
The regularity property (\ref{sol}) and the 
continuous dependence of solution on the initial data are shown by a usual way. 
This completes Theorem \ref{thm:lwp}. 
$\qed$

\subsection{Blowup criterion and scattering criterion} 

In this subsection we prove Theorems \ref{thm:sg}, \ref{prop:bc}, and \ref{prop:sc}.

\vskip2mm
\noindent
{\bf Proof of Theorem \ref{prop:bc}.} 
Assume for contradiction that
$T_{\mathrm{max}}<\I$
and $\| u \|_{S([0,T_{\mathrm{max}});\frac{\alpha-1}2)} <\I$.

{\bf Step 1}. We first show that the above assumption yields
$$\| u \|_{L([0,T_{\mathrm{max}}) ; \frac{\alpha-1}2)} <\I.$$
Fix $T$ so that $0< T < T_{\mathrm{max}}$.
Let $s_L(\alpha)$ be as in the previous section (see Remark \ref{rem:upperbound}).
If we take $\theta\in(0,1)$ so that $(\theta s_L(\alpha),\frac{\alpha-1}{2})$
is conjugate-acceptable then
it follows from
Proposition \ref{prop:ho_inho} that
\[
	\norm{u}_{L([0,T]; \frac{\alpha-1}2)}
	\le C\norm{ u_0 }_{\hat{L}^{\frac{\alpha-1}2}}
	+ C\norm{|u|^{\alpha -1} u }_{Y([0,T]; \theta s_L(\alpha),\frac{\alpha-1}2)}.
\]
Then, Lemma \ref{lem:nlest} (i) with $r=\frac{\alpha-1}2$
and Lemma \ref{lem:gHolder} give us
\[
	\norm{u}_{L([0,T]; \frac{\alpha-1}2)}
	\le C\norm{ u_0 }_{\hat{L}^{\frac{\alpha-1}2}}
	+ C \norm{u}_{S([0,T]; \frac{\alpha-1}2)}^{\alpha-\theta}
	\norm{u}_{L([0,T]; \frac{\alpha-1}2)}^\theta.
\]
By assumption, 
\[
	\norm{u}_{S([0,T];\frac{\alpha-1}2)}
	\le \norm{u}_{S([0,T_{\mathrm{max}}); \frac{\alpha-1}2)} <+\I
\]
for any $T \in (0, T_{\mathrm{max}})$.
Plugging this to the previous estimate, we see that 
there exist constants $A,B>0$ such that
\[
	\norm{u}_{L([0,T];\frac{\alpha-1}2)}
	\le A+ B\norm{u}_{L([0,T];\frac{\alpha-1}2)}^\theta
\]
for any $T \in (0, T_{\mathrm{max}})$, which gives us 
the desired bound since $\theta<1$.

{\bf Step 2}. Let $t_0 \in (0,T_{\mathrm{max}})$.
Since
\[
	u(t) = e^{-(t-t_0)\d_x^3}u(t_0)+\mu\int_{t_0}^t e^{-(t-t')\d_x^3} \d_x(|u|^{\alpha-1}u)(t') dt'
\]
for $t\in (0,T_{\mathrm{max}})$,
the above estimates yield the following bound on $e^{-(t-t_0)\d_x^3}u_0$:
\begin{eqnarray*}
\lefteqn{\norm{e^{-(t-t_0)\d_x^3}u(t_0)}_{S([t_0, T_{\mathrm{max}}) ;\frac{\alpha-1}2) \cap L([t_0,T_{\mathrm{max}};\frac{\alpha-1}2)}}\\
	&\le& 
	\norm{u}_{S([t_0, T_{\mathrm{max}}) ;\frac{\alpha-1}2) \cap L([t_0,T_{\mathrm{max}};\frac{\alpha-1}2)}
	\\
	& &+ C \norm{u}_{S([t_0, T_{\mathrm{max}});\frac{\alpha-1}2 )}^{\alpha-1}
	\norm{u}_{L([t_0,T_{\mathrm{max}});\frac{\alpha-1}2)}<\I.
\end{eqnarray*}

{\bf Step 3}. Let us now prove that we can extend the solution beyond $T_{\mathrm{max}}$.
Let $\delta$ be the constant given in Lemma \ref{lem:lwp}.
We see from the bound in the previous step  that there exists $t_0 \in (0,T_{\mathrm{max}}) $
 such that
\[
\norm{e^{-(t-t_0)\d_x^3}u(t_0)}_{S([t_0, T_{\mathrm{max}});\frac{\alpha-1}2 )}
	+ \norm{e^{-(t-t_0)\d_x^3}u(t_0)}_{L([t_0,T_{\mathrm{max}});\frac{\alpha-1}2)}
\le \frac{\delta}2.
\]
Hence, one can take $\tau>0$ so that
\[
	\norm{e^{-(t-t_0)\d_x^3}u(t_0)}_{S([t_0, T_{\mathrm{max}}+\tau);\frac{\alpha-1}2 )}
	+ \norm{e^{-(t-t_0)\d_x^3}u(t_0)}_{L([t_0,T_{\mathrm{max}}+\tau);\frac{\alpha-1}2)}
	\le {\delta}.
\]
Then, just as in the proof of Theorem \ref{thm:lwp} (or Lemma \ref{lem:lwp}),
we can construct a solution $u(t)$ to \eqref{gKdV}
in the interval $(-T_{\mathrm{min}},T_{\mathrm{max}}+\tau)$,
which contradicts to the definition of $T_{\mathrm{max}}$.
$\qed$

\vskip2mm
\noindent
{\bf Proof of Theorem \ref{prop:sc}.}
We first assume that $T_{\mathrm{max}}=+\I$ and 
$\| u \|_{S([0,\I);\frac{\alpha-1}2)} <\I$.
Then, as in the first step of the proof of Proposition \ref{prop:bc},
one obtains $\|u \|_{L([0,\I);\frac{\alpha-1}2)} <\I$.
Since $\{e^{-t\pt_{x}^{3}}\}_{t\in\rre}$ 
is isometry in $\hat{L}^{\frac{\alpha-1}{2}}$, 
it suffices to show that 
$\{e^{t\pt_{x}^{3}}u(t)\}_{t\in\rre}$ 
is a Cauchy sequence in $L^{\frac{\alpha-1}{2}}$ as $t\to\I$. 
Let $0<t_1<t_2$. By an argument similar to 
the proof of \eqref{eq:lwp_b2}, we obtain 
\begin{eqnarray*}
\left\|
e^{t_2\pt_{x}^{3}}u(t_2)
-e^{t_1\pt_{x}^{3}}u(t_1)
\right\|_{\hat{L}^{\frac{\alpha-1}{2}}
}&\le& 
C\||u|^{\alpha-1}u\|_{N([t_1,\I); \frac{\alpha-1}2)}\\
&\le& 
C\|u\|_{S([t_1,\I);\frac{\alpha-1}2)}^{\alpha-1}
\|u\|_{L([t_1,\I);\frac{\alpha-1}2)}
\\
& &\to0\quad as\ t_1\to\infty.
\end{eqnarray*}
Hence, we find that the solution to \eqref{gKdV}
scatters to a solution of the Airy equation 
as $t\to\infty$. 

Conversely, if $u(t)$ scatters forward in time then
we can choose $T>0$ so that
\[
	\norm{e^{-t\d_x^3} u_+}_{S([T,\I);\frac{\alpha-1}2)}
	+ \norm{e^{-t\d_x^3}u_+}_{L([T,\I);\frac{\alpha-1}2)}
	\le \frac{\delta}2,
\]
where $u_+ = \lim_{t\to\I} e^{t\d_x^3}u(t) \in \hat{L}^{\frac{\alpha-1}2}$ and
$\delta$ is the constant given in Lemma \ref{lem:lwp}.
Moreover, it holds for sufficiently large $t_0 \in [T,\I)$ that
\begin{align*}
	&\norm{e^{-t\d_x^{3}}(e^{t_0\d_x^3}u(t_0)- u_+)}_{S([T,\I);\frac{\alpha-1}2)}
	+ \norm{e^{-t\d_x^{3}}(e^{t_0\d_x^3}u(t_0)- u_+)}_{L([T,\I );\frac{\alpha-1}2)} \\
	&{}\le C \norm{e^{t_0\d_x^3}u(t_0)- u_+}_{\hat{L}^{\frac{\alpha-1}2}} \le \frac{\delta}2
\end{align*}
by means of \eqref{eq:mixed}.
We then see that
\[
	\norm{e^{-(t-t_0)\d_x} u(t_0)}_{S([T,\I);\frac{\alpha-1}2)}
	+ \norm{e^{-(t-t_0)\d_x}u(t_0)}_{L([T,\I );\frac{\alpha-1}2)}
	\le {\delta}.
\]
Then, Lemma \ref{lem:lwp} implies that
$\norm{u}_{S([T,\I);\frac{\alpha-1}2)}\le 2\delta$. $\qed$

\vskip2mm

\noindent
{\bf Proof of Theorem \ref{thm:sg}.} 
By \eqref{eq:mixed}, we have
\[
	\|e^{-t\d_x^3}u_0\|_{L(\R;\frac{\alpha-1}2)}
+\| e^{-t\d_x^3}u_0 \|_{S(\R;\frac{\alpha-1}2)} \le C \eps .
\]
Then, in light of Lemma \ref{lem:lwp}, we see that $u$ exists globally in time 
and satisfies $\norm{u}_S \le 2C \eps$, provided $\eps$ is small compared 
with the constant $\delta$ given in Lemma \ref{lem:lwp}. 
Proposition \ref{prop:sc} ensures that $u$ scatters for both time direction.
$\qed$

\subsection{Persistence of regularity}

In this subsection, we prove Theorems \ref{prop:reg1}, \ref{prop:reg2}, and then
\ref{thm:negativeE}.

\vskip2mm

\noindent
{\bf Proof of Theorem \ref{prop:reg1}.}
Let us prove that $u\in L(I;\frac{\alpha_0-1}2)$.
As in the proof of Lemma \ref{lem:lwp}, one deduces from Proposition \ref{prop:ho_inho} 
and Lemma \ref{lem:nlest} (i) that
\begin{align*}
	\norm{u}_{L(I;\frac{\alpha_0-1}2)}
	&{}\le C\norm{u_0}_{\hat{L}^{\frac{\alpha_0-1}2}}
	+ C\norm{|u|^{\alpha-1} u }_{N(I;\frac{\alpha_0-1}{2})} \\
	&{}\le C\norm{u_0}_{\hat{L}^{r_0}}
	+ C\norm{u}_{S(I;\frac{\alpha-1}2)}^{\alpha-1} \norm{u}_{L(I;\frac{\alpha_0-1}2)}.
\end{align*}
Since we already know $\norm{u}_{S(I;\frac{\alpha-1}2)} <\I$ by assumption,
we have the desired bound
\[
	\norm{u}_{L(I;\frac{\alpha_0-1}2)} \le 2C \norm{u_0}_{\hat{L}^{\frac{\alpha_0-1}2}}
\]
for sufficiently short interval $I$.
Then, again by Proposition \ref{prop:ho_inho},
\[
	\norm{u}_{L^\I_t (I;\hat{L}^{\frac{\alpha_0-1}2}_x) \cap X(I;s , \frac{\alpha_0-1}2)}
	\le C_s\norm{u_0}_{\hat{L}^{\frac{\alpha_0-1}2}}
	+ C_s\norm{u}_{S(I;\frac{\alpha-1}2)}^{\alpha-1} \norm{u}_{L(I;\frac{\alpha_0-1}2)} <+\I
\]
for any acceptable pair $(s,\frac{\alpha_0-1}2)$.
Finite time use of this argument yields the result.
$\qed$

\vskip2mm

\noindent
{\bf Proof of Theorem \ref{prop:reg2}.}
Let $0< \sigma < \alpha$.
Take a number $\eps$ so that $0< \eps < \min (1,\alpha - \sigma)$.
Since $|D_x|^\sigma$ commutes with $e^{-t\d_x^3}$ and since
$(\eps,2)$ is acceptable and conjugate-acceptable, we see
from Proposition \ref{prop:ho_inho} that
\begin{equation*}
	\norm{|D_x|^\sigma u(t)}_{X(I;\eps , 2)}\le
	C\norm{|D_x|^\sigma u_0}_{L^2} + C \norm{ |D_x|^\sigma (|u|^{\alpha-1} u )}_{Y(I;\eps,2)}.
\end{equation*}
Since $\sigma + \eps < \alpha$, arguing as in the proof of Lemma \ref{pp},
one sees that
\begin{eqnarray*}
\lefteqn{\norm{ |D_x|^\sigma (|u|^{\alpha-1} u )}_{Y(I;\eps,2)}}
\\
	&=& \norm{ |D_x|^{\sigma+\eps} (|u|^{\alpha-1} u)}_{
	L^{\tilde{p}(\eps ,2)}_xL^{\tilde{q}(\eps,2)}_t (I)} 
\\
	&\le& C\norm{u}_{L^{{p}(0,\frac{\alpha-1}2)}_xL^{{q}(0,\frac{\alpha-1}2)}_t (I)}^{\alpha-1} \norm{|D_x|^{\sigma+\eps} u}_{
	L^{{p}(\eps ,2)}_xL^{{q}(\eps,2)}_t (I)} 
	\\
	&=& C\norm{u}_{S(I;\frac{\alpha-1}2)}^{\alpha-1} \norm{|D_x|^{\sigma} u}_{
	X(I;\eps,2)}.
\end{eqnarray*}
Hence, we obtain an upper bound for $\norm{|D_x|^\sigma u}_{X(I;\eps,2)}$
for a small interval.
Then, the result follows as in Proposition \ref{prop:reg1}.

Next, let $-1 < \sigma < 0$. Set $\eps = -\sigma \in (0,1)$.
As in the previous case, we have
\begin{equation*}
	\norm{|D_x|^\sigma u(t)}_{X(I;\eps , 2)}\le
	C\norm{|D_x|^\sigma u_0}_{L^2} + C \norm{ |D_x|^\sigma (|u|^{\alpha-1} u )}_{Y(I;\eps,2)}
\end{equation*}
since $(\eps,2)$ is acceptable and conjugate-acceptable. Then,
\begin{align*}
	\norm{ |D_x|^\sigma (|u|^{\alpha-1} u )}_{Y(I;\eps,2)}
	&{}= \norm{|u|^{\alpha-1} u}_{L^{\tilde{p}(\eps ,2)}_xL^{\tilde{q}(\eps,2)}_t (I)} \\
	&{}\le \norm{u}_{S(I;\frac{\alpha-1}2)}^{\alpha-1} \norm{|D_x|^{\sigma} u}_{
	X(I;\eps,2)}
\end{align*}
by H\"older's inequality. The rest of the argument is the same.
$\qed$

\begin{remark}
In the above proposition, the upper bound $s<\alpha$ is natural
in view of the regularity which the nonlinearity $|u|^{\alpha-1} u$ possesses.
When $\alpha$ is an odd integer, that is, if $\alpha=5,7$, then the nonlinearity
$u^5$ or $u^7$ are analytic (in $u$) and so we can remove the upper bound
and treat all $s>0$. We omit the details. 
\end{remark}

\begin{remark}\label{rem:HsLWP}
By modifying the proof of Theorem \ref{prop:reg2},
we easily reproduce the local well-posedness in $\dot{H}^{s_\alpha}$ for $\alpha \ge5$.
More precisely, by Lemma \ref{lem:gHolder}, 
$$
	\norm{u}_{S(I;\frac{\alpha-1}{2})} \le 
	\norm{|D_x|^{s_\alpha} u}_{X(I;-1/4,2)}^{\frac{8}{5(\alpha-1)}}
	\norm{|D_x|^{\frac{2(9-\alpha)}{(5\alpha-13)(\alpha-1)}} u }_{L^{\frac{5\alpha-13}{2}}_{t,x}(I)}^{\frac{5\alpha -13}{5(\alpha-1)}}
$$
By Sobolev's embedding in space and Minkowski's inequality, 
\begin{align*}
	\norm{|D_x|^{\frac{2(9-\alpha)}{(5\alpha-13)(\alpha-1)}} u }_{L^{\frac{5\alpha-13}{2}}_{t,x}(I)}
	&{}\le C \norm{|D_x|^{s_\alpha - \frac{5\alpha-33}{4(5\alpha-13)}} u}_{L_t^{\frac{5\alpha-13}{2}} L_x^{\frac{4(5\alpha-13)}{5\alpha-17}}(I)}\\
	&{}\le C \norm{|D_x|^{s_\alpha} u}_{X(I;-\frac14 + \frac{5}{5\alpha-13},2)}
\end{align*}
Hence, estimating as in the proof of Theorem \ref{prop:reg2}, we obtain a closed estimate in
$|D_x|^{-s_\alpha} X(I;\eps , 2) \cap |D_x|^{-s_\alpha} 
X(I;-\frac14 + \frac{5}{5\alpha-13},2) \cap |D_x|^{-s_\alpha} X(I;-\frac14,2)$,
which yields local well-posedness in $\dot{H}^{s_\alpha}$.\footnote{
Strictly speaking, we should work with pairs
$(-\frac14 + \eta_1 ,2)$ and $(-\frac14 + \frac{5}{5\alpha-13} - \eta_2 ,2)$
for small $\eta_j = \eta_j(\alpha)>0$ because the critical case $q(-1/4,2) = \I$
 is excluded in Lemma \ref{lem:gHolder}. 
However, the modification is obvious.}
\end{remark}

\vskip2mm

We finally prove Theorem \ref{thm:negativeE}.

\vskip2mm

\noindent
{\bf Proof of Theorem \ref{thm:negativeE}.}
We suppose for contradiction that $u(t)$ scatters to $u_+ \in \hat{L}^{\frac{\alpha -1}{2}}$
as $t\to\I$.
Since $u_0 \in H^1$, Theorems \ref{prop:reg1} and \ref{prop:reg2} imply
that $u(t) \in C(\R; H^1)$.
Further, $u(t)$ scatters also in ${H}^1$ and so we see that
$\norm{\d_x u(t)}_{L^2} = \norm{\d_x e^{t\d_x^3} u(t)}_{L^2} \to \norm{u_+}_{\dot{H}^1}$
as $t\to\I$. 

On the other hand, by the Gagliardo-Nirenberg inequality and mass conservation, 
\[
	\norm{u(t)}_{L_{x}^{\alpha +1}}
	\le C \norm{u_0}_{L_{x}^2}^{\frac{2}{\alpha+1}} \norm{ |D_x|^{\frac{2}{3(\alpha-1)}} u(t)}_{L_x^{\frac{3(\alpha-1)}{2}}}^{\frac{\alpha -1}{\alpha +1}} .
\]
Since $u(t)$ scatters as $t\to\I$,
we see that $u\in X([0,\I); \frac2{3(\alpha-1)},\frac{\alpha-1}{2})$ as in the proof of Theorem \ref{prop:sc}.
Therefore, we can take a sequence $\{t_n\}_n$ with $t_n \to\I$ as $n\to\I$ so that
$\norm{u(t_n)}_{L^{\alpha +1}} \to 0$ as $n\to\I$.
Thus, by conservation of energy,
\[
	0 \ge E[u_0] =E[u(t_n)]= \frac12 \norm{\d_x u(t_n)}_{L^2}^2 -\frac{\mu}{\alpha + 1} \norm{u(t_n)}_{L^{\alpha+1}}^{\alpha+1} \to \frac12 \norm{u_+}_{\dot{H}^1}^2
\]
as $n\to \I$.
Hence, $E[u_0]<0$ yields a contradiction.
If $E[u_0]=0$ then we see that $u_+=0$, and so that
$\norm{u_0}_{L^2} = \norm{u_+}_{L^2} =0$.
This contradicts to $u_0 \neq 0$. $\qed$

\appendix


\section{Proof of Lemma \ref{ff}}

In this appendix we prove Lemma \ref{ff}. 
To prove this lemma, we need 
the following space-time bounds of 
the maximal function
\begin{eqnarray*}
({{\mathcal M}}u)(x)
=\sup_{R>0}\frac{1}{2R}
\int_{x-R}^{x+R}|u(y)|dy.
\end{eqnarray*}

\vskip2mm

\begin{lemma}\label{fs}
Let $I$ be an interval. Assume $1<p,q<\infty$.

\vskip1mm
\noindent
(i) There exists a positive 
constant $C$ depending on $p,q$ and $I$ such that 
\begin{eqnarray}
\|{{\mathcal M}}f\|_{L_{x}^{p}L_{t}^{q}(I)}
\le C\|f\|_{L_{x}^{p}L_{t}^{q}(I)} \label{fs1}
\end{eqnarray}
for any $f\in L_{x}^{p}L_{t}^{q}(I)$.

\vskip1mm
\noindent
(ii) There exists a positive 
constant $C$ depending on $p,q$ and $I$ such that 
\begin{eqnarray}
\|{{\mathcal M}}f_{k}\|_{L_{x}^{p}L_{t}^{q}\ell_{k}^{2}(I)}
\le C\|f_{k}\|_{L_{x}^{p}L_{t}^{q}\ell_{k}^{2}(I)} \label{fs2}
\end{eqnarray}
for any $\{f_{k}\}_{k}\in L_{x}^{p}L_{t}^{q}\ell_{k}^{2}(I)$.
\end{lemma}

\vskip2mm
\noindent
{\bf Proof of Lemma \ref{fs}.}
See \cite{FS} for (\ref{fs1}) and \cite[Lemma A.3 (e)]{KPV} for (\ref{fs2}). 
$\qed$

\vskip2mm
\noindent
{\bf Proof of Lemma \ref{ff}.} 
We follow \cite{Si} (see also \cite{RS}).
Let $\{\varphi_k(D_{x})\}_{k=-\I}^\I$ be a Littlewood-Paley decomposition 
with respect to $x$ variable. From \cite[Lemma A.3]{KPV}, we see 
\begin{equation}\label{eq:rep}
	\norm{|D_x|^s f}_{L^p_x L^q_t} 
	\sim \norm{2^{sk}\varphi_{k}(D_x)f}_{L^p_x L^q_t\ell_{k}^{2}}.
\end{equation}

{\bf Step 1.} 
Write $\mu = N+\beta$ with $N \in \Z$ and $\beta \in (0,1]$.
Remark that $N\ge1$ since $\mu>1$.
We first note that Taylor's expansion of $G$ gives us
\begin{align*}
	G(z)
	={}& \sum_{l=0}^{N-1}
	\frac{G^{(l)}(a)}{\ell !} (z-a)^l + \int_{a}^z \frac{(z-v)^{N-1}}{(N-1)!} G^{(N)}(v) dv \\
	={}& \sum_{l=0}^{N}
	\frac{G^{(l)}(a)}{\ell !} (z-a)^l + \int_{a}^z \frac{(z-v)^{N-1}}{(N-1)!} (G^{(N)}(v)-G^{(N)}(a)) dv\\
	={}& \sum_{l=0}^{N}\sum_{j=0}^l
	\frac{(-1)^{l-j} G^{(l)}(a) a^{l-j}}{(\ell-j) ! j!}
	z^{j} + \int_{a}^z \frac{(z-v)^{N-1}}{(N-1)!} (G^{(N)}(v)-G^{(N)}(a)) dv. \\
\end{align*}
Hence, applying the above expansion with $z=f(y)$ and $a=f(x)$,
\begin{equation}\label{eq:defT}
\begin{aligned}
	&\F^{-1}[\varphi_k \F G(f)](x) \\
	={}& c\int_{\R^n} (\F^{-1}\varphi_k) (x-y)G(f(y))dy \\
	={}& c\sum_{l=0}^{N}\sum_{j=0}^l 
	\frac{(-1)^{l-j} G^{(l)}(f(x)) (f(x))^{l-j} }{(\ell -j)! j!}
	\int_{\R^n} (\F^{-1}\varphi_k) (x-y) (f(y))^{j} dy\\
	&{} + c\int_{\R^n}  (\F^{-1}\varphi_k) (x-y)\int_{f(x)}^{f(y)} \frac{(f(y)-v)^{N-1}}{(N-1)!} 
	(G^{(N)}(v) -G^{(N)}(f(x)))dv  dy\\
	=:{}& T_{1,k} + T_{2,k}. 
\end{aligned}
\end{equation}
We first estimate $T_{1,k}$. Since $\int \F^{-1} \varphi_k(y) dy =\varphi_k(0) = 0$,
the summand in $T_{1,k}$ vanishes if $j=0$.
By the estimate
\[
	|G^{(l)}(f(x))| \le \norm{G}_{\mathrm{Lip}\mu} |f(x)|^{\mu-l},
\]
we have
\begin{align*}
	\norm{2^{sk}T_{1,k}}_{L^p_{x}L^q_{t}\ell^2_k}
	\le{}& C \norm{G}_{\mathrm{Lip}\mu} \sum_{j=1}^{N}
	\norm{ |f|^{\mu-j} \times 2^{sk} \varphi_k(D_{x})(f^j)| }_{L^p_{x}L^q_{t}\ell^2_k} \\
	\le{}& C \norm{G}_{\mathrm{Lip}\mu} \sum_{j=1}^{N} 
	\norm{f}_{L_{x}^{p_1}L_{t}^{q_{1}}}^{\mu-j}
	\norm{|D_{x}|^{s}(f^j)}_{L_{x}^{p_{2,j}}L^{q_{2,j}}_{t}},
\end{align*}
where 
\begin{eqnarray*}
\frac1p=\frac{\mu-j}{p_1}+\frac{1}{p_{2,j}},\qquad
\frac1q=\frac{\mu-j}{q_1}+\frac{1}{q_{2,j}}.
\end{eqnarray*}
Further, a recursive use of Lemma \ref{pp}
 yield
\[
	\norm{|D_{x}|^{s}(f^j)}_{L_{x}^{p_{2,j}}L_{t}^{q_{2,j}}}
	\le C_j \norm{f}_{L_{x}^{p_{1}}L_{t}^{q_{1}}}^{j-1} 
	\norm{|D_{x}|^{s}f}_{L_{x}^{p_{2}}L_{t}^{q_{2}}}
\]
for $j\ge2$, which completes the estimate of $T_{1,k}$.

Next, we estimate $T_{2,k}$. 
First note that
\begin{align*}
	\abs{\int_{f(x)}^{f(y)} \frac{(f(y)-v)^{N-1}}{(N-1)!} (G^{(N)}(v)-G^{(N)}(f(x))) dv}
	\le C \norm{G}_{\mathrm{Lip}\mu} |f(x)-f(y)|^{\mu}
\end{align*}
by definition of $\norm{G}_{\mathrm{Lip} \mu }$.
Further, for any $M>0$, there exists $C_M$ such that 
\[
	|(\F^{-1} \varphi_k) (x-y)|
	= 2^{k} |(\F^{-1} \varphi_0) (2^k(x-y))|
	\le C_M 2^{k} (1+ 2^k|x-y|)^{-M}.
\]
Therefore,
\begin{align*}
	|T_{2,k}| \le{}& C 2^{k} \norm{G}_{\mathrm{Lip}\mu}  \int_{\R^n}  \frac{|f(x)-f(y)|^{\mu}}{(1+2^k|x-y|)^M} dy \\
	\le{}& C \sum_{l=0}^\I 2^{k-l M} (I^\mu_{k-l} f)(x),
\end{align*}
where
\[
	I_k^\mu f(x) = \int_{|z|\le 2^{-k}} |f(x+z)-f(x)|^\mu dz.
\]
We now claim that
\begin{equation}\label{eq:claim}
\norm{ 2^{k(s+1)} (I^\mu_k f)}_{L^p_x L^q_t\ell_{k}^{2}}
\le C \norm{|D_{x}|^{s/\mu}f}_{L_{x}^{\mu p}L_{t}^{\mu q}}^\mu.
\end{equation}
This claim completes the proof. Indeed,
combining the above estimates, we see that
\begin{align*}
	\norm{2^{sk} T_{2,k}}_{L^p_xL^q_t\ell^2_k}
	\le C
	\sum_{l=0}^\I 2^{l(s-M+1)} \norm{ 2^{k(n+s)} (I^\mu_{k} f)}_{L^p_xL^q_t\ell_{k}^{2}}
	\le C \norm{|D_{x}|^{s/\mu}f}_{L_{x}^{\mu p}L_{t}^{\mu q}}^\mu,
\end{align*}
provided we choose $M>s+1$. By Lemma \ref{lem:gHolder}, we conclude that
\[
	\norm{|D_{x}|^{s/\mu}f}_{L_{x}^{\mu p}L_{t}^{\mu q}}^\mu
	\le 
	\norm{f}_{L_{x}^{p_{1}}L_{t}^{q_{1}}}^{1-\frac1\mu}
	\norm{|D_{x}|^{s}f}_{L_{x}^{p_{2}}L_{t}^{q_{2}}}^{\frac1\mu}.
\]

{\bf Step 2.}
We prove claim \eqref{eq:claim}.
Let $\Delta_h$ be a difference operator $\Delta_h f(x) = f(x+h)-f(x)$. 
Since $f = \sum_{m\in \Z} \varphi_{k+m}(D_{x}) f $ for any $k\in \Z$, one sees that
\begin{align*}
	\norm{ 2^{k(s+1)} (I^\mu_k f)(x)}_{L^p_{x}L^q_{t}\ell^2_k}
	={}& \norm{ 2^{ks} \int_{|z|\le 1} |\Delta_{2^{-k}z} f(x)|^\mu dz 
	}_{L^p_{x}L^q_{t}\ell^2_k} \\
	\le{}& \norm{ 2^{ks} \int_{|z|\le 1} |\Delta_{2^{-k}z} \sum_{m=-\I}^{-1} 
	\varphi_{k+m} (D)f(x)|^\mu dz }_{L^p_{x}L^q_{t}\ell^2_k} \\
	&{} + \norm{ 2^{ks} \int_{|z|\le 1} |\Delta_{2^{-k}z} \sum_{m=0}^{\I} 
	\varphi_{k+m} (D)f(x)|^\mu dz }_{L^p_{x}L^q_{t}\ell^2_k}\\
	=:{}& A + B.
\end{align*}

We estimate $A$.
Take $a\in (1/\mu,1)$.
Let $k\in \Z$.
If $m<0$ and $|h| \le 2^{-k}$ then
we have
\begin{align*}
	|\Delta_h \F^{-1} [\varphi_{k+m} \F f](x)|
	\le{}& |h| | \nabla( \F^{-1} [\varphi_{k+m} \F f])(x+\theta h)|\\
	\le{}& 2^{m} \sup_{|y|\le 2^{-k}} |
	(\nabla \F^{-1}[ {\varphi}_{0}\F[f(\frac{\cdot}{2^{k+m}})]])(2^{k+m}(x-y))| \\
	\le{}& C_a 2^{m} \sup_{ y\in \R} \frac{(\nabla \F^{-1}[ {\varphi}_{0}
	\F[f(\frac{\cdot}{2^{k+m}})]])(2^{k+m}(x-y))|}{1+|2^{k+m}y|^a } \\
	\le{}& C_a 2^{m} \sup_{y\in \R} \frac{| \F^{-1} [ {\varphi}_{k+m} \F f](x-y)|}{1+|2^{k+m}y|^a}
\end{align*}
for any  $x\in \R$, where we have used the estimate
\[
	\sup_{ y\in \R} \frac{|\nabla \F^{-1} [ {\varphi}_{0} \F f](x-y)|}{1+|y|^a }
	\le 
	C  \sup_{y\in \R} \frac{| \F^{-1} [ {\varphi}_{0} \F f](x-y)|}{1+|y|^a}
\]
(see \cite[Proposition 2.1.6/2 (i)]{RS}) to obtain the last line.
We define the Peetre-Fefferman-Stein maximal function by 
\[
	\varphi^{*,a}_{j} f (x) := \sup_{y\in \R} 
	\frac{| \F^{-1} [ {\varphi}_{j} \F f](x-y)|}{1+|2^{j}y|^a}.
\]
By the above estimates, we have
\begin{align*}
	A \le {}&
	C\norm{ 2^{ks}  \sum_{m=-\I}^{-1} \sup_{|z|\le 1} |\Delta_{2^{-k}z} \varphi_{k+m} (D)f(x)|^\mu
	}_{L^p_{x}L^q_{t}\ell^2_k} \\
	\le{}& C  \sum_{m=-\I}^{-1} 2^{m\mu}
	\norm{ 2^{k\frac{s}\mu}\varphi^{*,a}_{k+m} f  }_{L^{\mu p}_{x}L^{\mu q}_{t}\ell^{2\mu}_k}^{\mu} \\
	\le{}& C  \sum_{m=-\I}^{-1} 2^{m(\mu-s)}
	\norm{ 2^{(k+m)\frac{s}\mu}\varphi^{*,a}_{k+m} f  }_{L^{\mu p}_{x}L^{\mu q}_{t}\ell^{2\mu}_k}^{\mu}\\
	\le{}& C 
	\norm{ 2^{k\frac{s}\mu}\varphi^{*,a}_{k} f  }_{L^{\mu p}_{x}L^{\mu q}_{t}\ell^{2\mu}_k}^{\mu},
\end{align*}
where we used the fact that $s<\mu$. 
Since $(\varphi_{k}^{*,a}f)(x)=(\varphi_{0}^{*,a}(\tilde{\varphi}_{k}(D_{x})f)
(\frac{\cdot}{2^{k}}))(2^{k}x)$, \cite[Lemma 2.3.6]{Tri} yields
\begin{eqnarray*}
(\varphi_{k}^{*,a}f)(x)
\le C\{{{\mathcal M}}[(\tilde{\varphi}_{k}(D_{x})f)^{\frac1a}]\}^{a}(x),
\end{eqnarray*}
where $\tilde\varphi_k = \sum_{i=-1}^1 \varphi_{k+i}$.
Since $1/\mu<a<1$,   
\eqref{fs2}, the embedding $\ell^2 \hookrightarrow \ell^q$ ($2<q\le \I$), and \eqref{eq:rep}
lead us to
\begin{eqnarray*}
	\norm{ 2^{k\frac{s}\mu}\varphi^{*,a}_{k} f   
	}_{L^{\mu p}_{x}L^{\mu q}_{t}\ell^{2\mu}_k}
	&\le& C\norm{2^{k\frac{s}{a\mu}}{{\mathcal M}}[(\tilde{\varphi}_{k}(D_{x})f)^{\frac1a}]
	}_{L^{a\mu p}_{x}L^{a\mu q}_{t}\ell^{2a\mu}_k}^{a}\\
	&\le& C\norm{2^{k\frac{s}{a\mu}}(\tilde{\varphi}_{k}(D_{x})f)^{\frac1a}
	}_{L^{a\mu p}_{x}L^{a\mu q}_{t}\ell^{2}_k}^{a}\\
	&\le& C\norm{2^{k\frac{s}{\mu}}\tilde{\varphi}_{k}(D_{x})f
	}_{L^{\mu p}_{x}L^{\mu q}_{t}\ell^{\frac{2}{a}}_k}\\
	&\le& C \norm{|D_{x}|^{s/\mu}f}_{L_{x}^{\mu p}L_{t}^{\mu q}}.
\end{eqnarray*}

Let us proceed to the estimate of $B$.
We first note that
\begin{align*}
	&\int_{|z|\le 1} \left|\Delta_{2^{-k}z} \sum_{m=0}^{\I} \varphi_{k+m} (D)f(x)\right|^\mu dz \\
	&{}= \int_{|z|\le 1} \left| \sum_{m=0}^{\I} 2^{-\frac{\eps}{\mu} m} 2^{\frac{\eps}{\mu} m}\Delta_{2^{-k}z} \varphi_{k+m} (D)f(x)\right|^\mu dz \\
	&{}\le C_{\eps}\int_{|z|\le 1} \sum_{m=0}^{\I} 2^{{\eps} m}|\Delta_{2^{-k}z} \varphi_{k+m} (D)f(x)|^\mu dz  \\
	&{}= C_\eps \sum_{m=0}^{\I} 2^{{\eps} m} \int_{|z|\le 1}|\Delta_{2^{-k}z} \varphi_{k+m} (D)f(x)|^\mu dz \\
	&{}\le C \sum_{m=0}^{\I} 2^{{\eps} m} \(\sup_{|z|\le 1}|\Delta_{2^{-k}z} \varphi_{k+m} (D)f(x)|\)^{\mu(1-\lambda)} \\
&{}\quad \times \int_{|z|\le 1}|\Delta_{2^{-k}z} \varphi_{k+m} (D)f(x)|^{\mu\lambda} dz,
\end{align*}
where $\lambda \in (0,1)$.
For $m\ge 0$ and $|h| \le 2^{-k}$, the triangle inequality gives us
\begin{align*}
	|\Delta_h \F^{-1} [\varphi_{k+m} \F f](x)|
	&{}\le 2 \sup_{|y|\le 2^{-k}} |\F^{-1} [\varphi_{k+m} \F f](x-y)| \\
	&{}\le C 2^{ma}\varphi^{*,a}_{k+m} f (x),
\end{align*}
where $a\in (1/\mu,1)$. Further,
\[
	\int_{|z|\le 1}|\Delta_{2^{-k}z} \varphi_{k+m} (D_{x})f(x)|^{\mu\lambda} dz
	\le C \mathcal{M}[|\varphi_{k+m} (D)_{x}f|^{\mu\lambda}](x).
\]
Plugging these inequality, one deduces 
from H\"older's inequality,
the embedding $\ell^2 \hookrightarrow \ell^q$ ($2<q\le \I$),
\eqref{fs2}, and \eqref{eq:rep} 
that
\begin{align*}
	B \le{}& C \norm{2^{sk} \sum_{m=0}^\I 2^{m\eps} \mathcal{M}[|\varphi_{k+m} (D_{x})f|^{\mu\lambda}]
	2^{ma\mu(1-\lambda)}(\varphi^{*,a}_{k+m} f )^{\mu(1-\lambda)}
	}_{L^p_{x}L^q_{t}\ell^2_k} \\
	\le{}& C \sum_{m=0}^\I 2^{m(\eps+a\mu(1-\lambda))}\norm{2^{sk} \mathcal{M}[|\varphi_{k+m} (D_{x})f|^{\mu\lambda}]
	(\varphi^{*,a}_{k+m} f )^{\mu(1-\lambda)}
	}_{L^p_{x}L^q_{t}\ell^2_k} \\
	\le{}& C \sum_{m=0}^\I 2^{m(\eps+a\mu(1-\lambda) -s)}\norm{ \mathcal{M}[|2^{\frac{s}{\mu}k}\varphi_{k} (D_{x})f|^{\mu\lambda}] }_{L^{\frac{p}{\lambda}}_xL^{\frac{q}{\lambda}}_t\ell^{\frac{2}{\lambda}}_k}
	\norm{2^{\frac{s}{\mu}k}
	\varphi^{*,a}_{k} f}_{L^{\mu p}_xL^{\mu q}_t\ell^{2\mu}_k}^{\mu(1-\lambda)}\\
	\le{}&
	C \sum_{m=0}^\I 2^{m(\eps+a\mu(1-\lambda) -s)}
	\norm{|D_{x}|^{s/\mu}f}_{L_{x}^{\mu p}L_{t}^{\mu q}}^\mu \\
	\le{}&\norm{|D_{x}|^{s/\mu}f}_{L_{x}^{\mu p}L_{t}^{\mu q}}^\mu.
\end{align*}
as long as $\eps+a\mu(1-\lambda) -s<0$. 
Since $a\in(1/\mu,1)$,
we are able to choose $\lambda \in (0,1)$ and $\eps>0$ suitably.
Thus, the proof is completed. $\qed$


\section{Inclusion relations of $\hat{L}^{r}$}\label{sec:embedding}

In this appendix we briefly summarize some inclusion relations between $\hat{L}^r$
and other frequently used spaces such as Lebesgue space or Sobolev space.
Here, $\dot{H}^{0,s}=\dot{H}^{0,s}(\R)$ stands for a weighted $L^2$ space with norm $\norm{f}_{\dot{H}^{0,s}}= \||x|^s f \|_{L^2}$.

\begin{lemma}\label{lem:basic_inclusion} 
We have the following.

\vskip1mm
\noindent
(i) $ {L}^r\hookrightarrow\hat{L}^r$ if $1\le r \le 2$ and
$\hat{L}^r\hookrightarrow L^r$ if $2\le r \le \I$.

\vskip1mm
\noindent
(ii) $\dot{H}^{0,\frac1r-\frac12}\hookrightarrow\hat{L}^r$ 
if $1<r \le 2$
and $\hat{L}^r \hookrightarrow \dot{H}^{0,\frac1r-\frac12} $ if $2 \le r <\I$.

\vskip1mm
\noindent
(iii) $\hat{L}^r \hookrightarrow \dot{B}^{\frac12-\frac1r}_{2,r'}$ 
if $1\le r \le 2$
and $\dot{B}^{\frac12-\frac1r}_{2,r'} \hookrightarrow \hat{L}^r$ if $2 \le r \le \I$
\end{lemma}

\vskip2mm
\noindent
{\bf Proof of Lemma \ref{lem:basic_inclusion}.}
The first assertion follows from the Hausdorff-Young inequality.
The Sobolev embedding (in Fourier side) yields the second.
We omit details.

The third is also immediate from the H\"older inequality.
Indeed, if $2\le r \le \I$ then
\[
	\norm{\hat{f}}_{L^{r'}(\{2^{n} \le |\xi| \le 2^{n+1}\})} \le
	C 2^{n(\frac12 -\frac1r)} \norm{\hat{f}}_{L^{2}(\{2^{n} \le |\xi| \le 2^{n+1}\})}
\]
for any $n\in \Z$. Taking $\ell^{r'}_n$ norm, we obtain the desired embedding.
The case $1\le r \le 2$ follows in the same way.
$\qed$

\vskip3mm
Let $\dot{H}^s=\dot{H}^s(\R)$ be a homogeneous Sobolev space
with norm $\norm{f}_{\dot{H}^s}= \||\xi|^s \hat{f} \|_{L^2}$.
Notice that the above inclusion is the same as for $\dot{H}^{\frac12-\frac1r}$.
Namely, we can replace $\hat{L}^r$ with $\dot{H}^{\frac12-\frac1r}$ in Lemma \ref{lem:basic_inclusion}.
However, there is no inclusion between $\hat{L}^r$ and $\dot{H}^{\frac12-\frac1r}$ for $r\neq2$.
\begin{lemma}\label{lem:ws}
For $1\le r\le \I$ ($r\neq 2$), $\hat{L}^r \not\hookrightarrow \dot{H}^{\frac12-\frac1r}$
and $\dot{H}^{\frac12-\frac1r} \not\hookrightarrow \hat{L}^r $.
\end{lemma}

\vskip2mm
\noindent
{\bf Proof of Lemma \ref{lem:ws}.}
If $2<r\le \I$, we have the following counter examples:
Let us define $f_n(x)$ by
$\hat{f_n}(\xi)=1$ for $n \le \xi \le n+1$ and $\hat{f_n}(\xi) =0$ elsewhere.
Then, $f_n(x)$ satisfies
$\norm{f_n}_{\dot{H}^{\frac12-\frac1r}} \to \I$ as $n\to\I$,
while $\norm{f_n}_{\hat{L}^r} = 1$.
Hence. $\hat{L}^r \not\hookrightarrow \dot{H}^{\frac12-\frac1r}$.
On the other hand, for some $p\in(1/2,1/r')$, 
take $g_n(x)$ ($n\ge3$) so that
$\hat{g_n}(\xi)=\xi^{-1/r'}|\log \xi|^{-p}$ for $1/n \le \xi \le 1/2$ and $\hat{g_n}(\xi) =0$ elsewhere.
Then, $\norm{g_n}_{\dot{H}^{\frac12-\frac1r}}$ is bounded
but $\norm{g_n}_{\hat{L}^r} \to \I$ as $n\to\I$.
This shows $\dot{H}^{\frac12-\frac1r} \not\hookrightarrow \hat{L}^r $.

The case $1<r <2$ follows by duality.

Let us consider the case $r=1$.
We note that $\delta_0(x) \in \hat{L}^1 \setminus \dot{H}^{-\frac12}$, where $\delta_0(x)$
is the Dirac delta function.
Therefore, $\hat{L}^1 \not\hookrightarrow \dot{H}^{-\frac12}$.
On the other hand, $f_n (x) = (\log (1+1/n))^{-1} \F^{-1} [{\bf 1}_{\{1\le \xi \le 1+1/n\}}] (x)$
is a counter example for $\dot{H}^{-\frac12} \not\hookrightarrow \hat{L}^1 $.
$\qed$

\vskip3mm
\noindent {\bf Acknowledgments.}
The authors express their deep gratitude to Professor Yoshio Tsutsumi for
valuable comments on preliminary version of the manuscript.
They would particularly like to thank Professor Gr\"{u}nrock for drawing their attention to the article \cite{G1} in which a more general version of
Lemma \ref{ST} is proved. 
The part of this work was done while the authors were 
visiting at Department of Mathematics at the University of
California, Santa Barbara whose hospitality they gratefully
acknowledge.
S.M. is partially supported by JSPS, Grant-in-Aid for Young Scientists (B) 24740108.
J.S. is partially supported by JSPS, Strategic Young Researcher Overseas
Visits Program for Accelerating Brain Circulation and by MEXT,
Grant-in-Aid for Young Scientists (A) 25707004.



\begin{thebibliography}{30}

\bibitem{B} Bourgain J., \textit{
Fourier restriction phenomena for certain lattice subsets and
applications to nonlinear evolution equations II. The KdV
equation}. Geom, Funct. Anal. {\bf 3} (1993), 209--262.

\bibitem{BK} Buckmaster T. and Koch H., 
\textit{The Korteweg-de Vries equation at $H^{-1}$ regularity}, 
Ann. Inst. H. Poincar\'e Anal. Non Lin\'eaire (C), To appear. 


\bibitem{CK} Christ M. and Kiselev A., 
\textit{Maximal functions associated to filtrations},
J. Funct. Anal. {\bf179} (2001) 409--425. 

\bibitem{CW} 
Christ F.M. and Weinstein M.I., 
\textit{Dispersion of small amplitude solutions of 
the generalized Korteweg-de Vries equation}. 
J. Funct. Anal. {\bf 100} (1991) 87--109.

\bibitem{F}
Fefferman C., 
\textit{Inequalities for strongly singular 
convolution operators}. 
Acta Math. {\bf 124} (1970) 9--36. 

\bibitem{FS}
Fefferman C. and Stein E.M., 
\textit{Some maximal inequalities}. 
Amer. J. Math. {\bf 93} (1971) 107--115.

\bibitem{Fz}
Fernandez D.L.,
\textit{Vector-valued singular integral operators 
on $L^p$-spaces with
mixed norms and applications}.
Pacific J. Math.  {\bf 129} (1987),  no. 2, 257--275.

\bibitem{G}
Gr\"{u}nrock A., \textit{
A bilinear Airy estimate with application to gKdV-3}, 
Differential Integral Equations,
{\bf 18} (2005), 1333--1339.

\bibitem{G1}
Gr\"{u}nrock A., 
\textit{An improved local well-posedness result 
for the modified KdV equation}. 
Int. Math. Res. Not. {\bf 2004} (2004), 
3287--3308.

\bibitem{G2}
Gr\"{u}nrock A., 
\textit{Bi- and trilinear Schr\"{o}dinger estimates 
in one space dimension with applications to cubic NLS 
and DNLS. Int}. Math. Res. Not. {\bf 2005}, (2005), 
2525--2558. 

\bibitem{GV}
Gr\"{u}nrock A. and Vega L., 
\textit{Local well-posedness for the modified KdV equation 
in almost critical $\hat{H}^{r}$-spaces}. 
Trans. Amer. Math. Soc. {\bf 361} (2009), 5681--5694.

\bibitem{Guo} Guo Z., 
\textit{Global well-posedness of Korteweg-de Vries 
equation in $H^{-3/4}(\rre)$}. 
J. Math. Pures Appl. (9), {\bf 91} (2009), 583--597.




\bibitem{H}
Hirschman I., Jr.,
\textit{A convexity theorem for certain groups 
of transformations}. J. Analyse Math. {\bf 2} 
(1953), 209--218.

\bibitem{HT} 
Hyakuna R. and Tsutsumi M., 
\textit{On existence of global solutions of 
Schr\"{o}dinger equations with subcritical 
nonlinearity for $\hat{L}^{p}$-initial data}. 
Proc. Amer. Math. Soc. {\bf 140} (2012), 
3905--3920.


\bibitem{Kato} 
Kato T., \textit{On the Cauchy problem 
for the (generalized) KdV equation}.
Advances in Math. Supplementary studies, 
Studies in Applied Mathematics {\bf 8} 
(1983), 93--128.

\bibitem{KPV1}
Kenig C.E., Ponce G. and Vega L., \textit{
Oscillatory integrals and regularity
of dispersive equations}. Indiana Univ.math J. 
{\bf 40} (1991), 33--69.

\bibitem{KPV}
Kenig C.E., Ponce G. and Vega L., \textit{
Well-posedness and scattering results for the generalized 
Korteweg-de Vries equation via the contraction principle}. 
Comm. Pure Appl. Math. {\bf 46} (1993) 527--620.

\bibitem{KPVq} 
Kenig C.E., Ponce G. and Vega L., \textit{
 A bilinear estimate with applications
to the KdV equation}. J. Amer. Math. Soc. 
{\bf 9} (1996), 573--603.

\bibitem{KeRu} 
Kenig C.E. and Ruiz A., 
\textit{A strong type $(2,2)$ 
estimate for the maximal
function associated to the 
Schr\"{o}dinger equation}.
Trans. Amer. Math. Soc. {\bf 280}
(1983), 239--246. 


\bibitem{KN}
Kishimoto N.,\textit{
Well-posedness of the Cauchy problem for 
the Korteweg-de Vries equation at the 
critical regularity}. 
Differential Integral Equations {\bf 22} (2009), 
447--464. 

\bibitem{KM}  Koch H. and Marzuola J.L., 
\textit{Small data scattering and soliton stability in 
$\dot{H}^{-\frac16}$ 
for the quartic KdV equation}. 
Anal. PDE {\bf 5} (2012), 145--198.

\bibitem{KV}
Korteweg D. J. and de Vries G., 
\textit{On the change of form of long waves advancing in a rectangular canal, 
and on a new type of long stationary waves}, Philos. Mag. {\bf 39} (1895), 422-443.

\bibitem{L} Lamb G.L.Jr., \textit{
Solitons on moving space curves}. 
J. Math. Phys. {\bf 18} (1977), 1654--1661. 

\bibitem{MR1} Molinet L. and Ribaud F., \textit{
On the Cauchy problem for the generalized Korteweg-de Vries 
equation}. Comm. Partial Differential Equations {\bf 28} (2003), 
2065--2091.

\bibitem{MR} Molinet L. and Ribaud F., \textit{ 
Well-posedness results for the generalized 
Benjamin-Ono equation with small initial data}. 
J. Math. Pures Appl. {\bf83} (2004), 277--311.



\bibitem{RS} Runst T. and Sickel W., 
\textit{``Sobolev spaces of fractional order, Nemytskij operators, and nonlinear partial 
differential equations''}. Walter de Gruyter \& Co., Berlin, (1996).



\bibitem{Si}
Sickel W., 
\textit{On boundedness of superposition operators in spaces of Triebel-Lizorkin type}.
Czechoslovak Math. J. {\bf 39(114)} (1989), 323--347. 

\bibitem{S}
Stein E. M.,
\textit{Interpolation of linear operators}.
Trans. Amer. Math. Soc. {\bf 83}  (1956), 482--492.


\bibitem{Strunk}
Strunk N., \textit{
Well-posedness for the supercritical 
gKdV equation}. 
Commun. Pure Appl. Anal. {\bf 13} 
(2014), 527--542.

\bibitem{T} 
Tao T., \textit
{Scattering for the quartic generalised Korteweg-de 
Vries equation}. 
J. Differential Equations {\bf 232} (2007), 623--651.

\bibitem{TVV}
Tao T., Vargas A., and Vega L.,
\textit{A bilinear approach to the restriction and Kakeya conjectures}.
J. Amer. Math. Soc. {\bf 11} (1998), 967--1000.

\bibitem{To}
Tomas P. A., 
\textit{A restriction theorem for the Fourier transform}.
Bull. Amer. Math. Soc. {\bf 81} (1975), 477--478.

\bibitem{Tri} Triebel H., \textit{``Theory of function spaces''}.
Monographs in Mathematics, {\bf 78} 
Birkh\"auser Verlag, Basel (1983).

\end{thebibliography}
\end{document}